\documentclass[leqno,11pt]{article}

\usepackage{tikz}
\usepackage{tikz-cd}
\usepackage{amsmath}

\usepackage[utf8]{inputenc}
\usepackage[T1]{fontenc}
\usepackage{microtype}

\usepackage{calc}
\usepackage[a4paper]{geometry}

\usepackage{amsmath}
\usepackage{amsthm}
\usepackage{tikz-cd}
\usetikzlibrary{arrows} 
\tikzset{
  commutative diagrams/.cd, 
  arrow style=tikz, 
  diagrams={>=stealth}
}
\usetikzlibrary{matrix,decorations.pathreplacing,calc}

\usepackage{textcomp}
\usepackage[sb]{libertine}
\usepackage[varqu,varl]{zi4}%
\usepackage[libertine,bigdelims,vvarbb]{newtxmath}

\usepackage[scr=boondox,cal=euler]{mathalfa}

\usepackage{marvosym}

\usepackage{slashed}
\usepackage{esint} 
\usepackage[english]{babel}

\usepackage{imakeidx}
\indexsetup{noclearpage}
\makeindex[intoc]

\usepackage{csquotes}
\usepackage[
  backend=biber,
  hyperref=true,
  backref=true,
  isbn=false,
  doi=true,
  natbib=true,
  eprint=true,
  useprefix=true,
  maxcitenames=99,
  maxbibnames=99,  
  maxalphanames=99, 
  minalphanames=99,
  safeinputenc,
  style=alphabetic,
  citestyle=alphabetic,
  block=space,
  datamodel=preamble/ext-eprint,
  sorting=nyt
]{biblatex}
\definecolor{darkblue}{HTML}{0000A6}

\usepackage[
  bookmarksnumbered = true,
  hypertexnames     = false,
  colorlinks        = true,
  citecolor         = darkblue,
  linkcolor         = darkblue,
  urlcolor          = darkblue,
  breaklinks
]{hyperref}

\DeclareFieldFormat{url}{%
  \href{#1}{\ComputerMouse}
}
\DeclareFieldFormat{doi}{%
  \mkbibacro{DOI}\addcolon\space\href{https://doi.org/#1}{#1}
}
\makeatletter
\DeclareFieldFormat{arxiv}{%
  arXiv\addcolon\space\href{http://arxiv.org/\abx@arxivpath/#1}{#1}
}
\makeatother
\DeclareFieldFormat{mr}{%
  MR\addcolon\space\href{http://www.ams.org/mathscinet-getitem?mr=MR#1}{#1}
}
\DeclareFieldFormat{zbl}{%
  Zbl\addcolon\space\href{http://zbmath.org/?q=an:#1}{#1}
}
\renewbibmacro*{eprint}{%
  \printfield{arxiv}%
  \newunit\newblock
  \printfield{mr}%
  \newunit\newblock
  \printfield{zbl}%
  \newunit\newblock
  \iffieldundef{eprinttype}
  {\printfield{eprint}}
  {\printfield[eprint:\strfield{eprinttype}]{eprint}}
}

\AtEveryBibitem{%
  \clearlist{address}%
}
\DeclareFieldFormat[article,inproceedings,inbook,incollection,thesis]{title}{\textit{#1}}
\renewbibmacro{in:}{}
\newcommand{\printreferences}{\raggedright\printbibliography[heading=bibintoc]}
\usepackage[inline,shortlabels]{enumitem}

\usepackage{subcaption}

\usepackage{etoolbox}
\ifundef{\abstract}{}{\patchcmd{\abstract}%
    {\quotation}{\quotation\noindent\ignorespaces}{}{}}

\usepackage[super]{nth}

\usepackage{thmtools}

\numberwithin{equation}{section}
\renewcommand{\qedsymbol}{$\blacksquare$}

\newcommand{\CorollaryQED}{\qedsymbol}

\newcommand{\ConjectureQED}{$\square$}
\newcommand{\SituationQED}{$\times$}
\newcommand{\DefinitionQED}{$\spadesuit$}
\newcommand{\NotationQED}{$\blacktriangleright$}
\newcommand{\ExampleQED}{$\bullet$}
\newcommand{\RemarkQED}{$\clubsuit$}

\declaretheorem[numberlike=equation,]{theorem}
\declaretheorem[numbered=no,name=Theorem]{theorem*}
\declaretheorem[numberlike=equation,name=Lemma]{lemma}
\declaretheorem[numberlike=equation,name=Proposition]{prop}
\declaretheorem[numberlike=equation,name=Corollary,qed=\CorollaryQED]{cor}

\declaretheorem[numberlike=equation,name=Definition,style=definition,qed=\DefinitionQED]{definition}
\declaretheorem[numbered=no,name=Definition,style=definition,qed=\DefinitionQED]{definition*}
\declaretheorem[numberlike=equation,name=Notation,style=definition,qed=\NotationQED]{notation}

\declaretheorem[numberlike=equation,style=definition,qed=\ExampleQED]{example}

\declaretheorem[numberlike=equation,style=remark,qed=\RemarkQED]{remark}
\declaretheorem[numbered=no,style=remark,name=Remark,qed=\RemarkQED]{remark*}

\def\makeautorefname#1#2{\AtBeginDocument{\expandafter\def\csname#1autorefname\endcsname{#2}}}
\makeautorefname{table}{Table}        
\makeautorefname{chapter}{Chapter}
\makeautorefname{section}{Section}
\makeautorefname{subsection}{Section}
\makeautorefname{subsubsection}{Section}
\makeautorefname{footnote}{Footnote}
\AtBeginDocument{\def\itemautorefname~#1\null{(#1)\null}}
\AtBeginDocument{\def\equationautorefname~#1\null{(#1)\null}}

\numberwithin{substep}{step}
\makeautorefname{step}{Step}
\makeautorefname{substep}{Step}

\makeautorefname{case}{Case}
\makeautorefname{substep}{Step}

\setlist[description]{leftmargin=!,labelindent=1em}
\setlist[enumerate]{label={\rm (\arabic*)},ref=\arabic*}
\setlist[enumerate,2]{label={\rm (\alph*)},ref=\theenumi.\alph*}
\setlist[enumerate,3]{label={\rm (\roman*)},ref=\theenumii.\roman*}

\let\C\undefined

\usepackage{bm}
\usepackage{mathtools} 
\usepackage{stmaryrd} 

\DeclareFontFamily{U}{mathx}{\hyphenchar\font45}
\DeclareFontShape{U}{mathx}{m}{n}{
      <5> <6> <7> <8> <9> <10>
      <10.95> <12> <14.4> <17.28> <20.74> <24.88>
      mathx10
      }{}
\DeclareSymbolFont{mathx}{U}{mathx}{m}{n}
\DeclareFontSubstitution{U}{mathx}{m}{n}
\DeclareMathAccent{\widecheck}{0}{mathx}{"71}
\DeclareMathAccent{\wideparen}{0}{mathx}{"75}

\DeclareMathOperator{\Diff}{Diff}
\DeclareMathOperator{\End}{End}

\DeclareMathOperator{\graph}{graph}
\DeclareMathOperator{\HF}{\HF}

\DeclareMathOperator{\Hom}{Hom}

\DeclareMathOperator{\coker}{coker}

\DeclareMathOperator{\im}{im}
\DeclareMathOperator{\ind}{index}

\DeclarePairedDelimiter{\norm}{\|}{\|}

\DeclarePairedDelimiterX{\inp}[2]{\langle}{\rangle}{#1, #2}

\DeclarePairedDelimiter{\abs}{\lvert}{\rvert}

\def\({\left(}
\def\){\right)}
\def\<{\left\langle}
\def\>{\right\rangle}

\newcommand{\C}{{\mathbf{C}}}

\newcommand{\N}{{\mathbf{N}}}

\newcommand{\R}{\mathbf{R}}

\newcommand{\SU}{\mathrm{SU}}

\newcommand{\Vect}{\mathrm{Vect}}

\newcommand{\Z}{\mathbf{Z}}

\newcommand{\hol}{\mathrm{hol}}

\newcommand{\id}{\mathbf 1}

\newcommand{\loc}{\mathrm{loc}}

\newcommand{\ob}{\mathrm{ob}}

\newcommand{\qandq}{\quad\text{and}\quad}

\renewcommand{\epsilon}{\varepsilon}

\newcommand{\vol}{\mathrm{vol}}

\renewcommand{\Re}{\operatorname{Re}}

\renewcommand{\emptyset}{\varnothing}
\renewcommand{\setminus}{{\backslash}}

\renewcommand{\leq}{\leqslant}
\renewcommand{\geq}{\geqslant}

\newcommand{\tn}{\otimes}

\makeatletter
\renewcommand*\env@matrix[1][*\c@MaxMatrixCols c]{%
  \hskip -\arraycolsep
  \let\@ifnextchar\new@ifnextchar
  \array{#1}}

\renewcommand\xleftrightarrow[2][]{%
  \ext@arrow 9999{\longleftrightarrowfill@}{#1}{#2}}
\newcommand\longleftrightarrowfill@{%
  \arrowfill@\leftarrow\relbar\rightarrow}
\makeatother

\newcommand{\bD}{{\mathbf{D}}}

\newcommand{\cD}{\mathcal{D}}

\newcommand{\cL}{\mathcal{L}}
\newcommand{\cM}{\mathcal{M}}

\newcommand{\fL}{{\mathfrak L}}

\author{Gorapada Bera}
\title{Deformations of asymptotically cylindrical associative submanifolds}
\date{\vspace{-5ex}}
\begin{document}
\maketitle
\begin{abstract}
This article develops the deformation theory of asymptotically cylindrical (ACyl) associative submanifolds in ACyl $G_2$-manifolds, laying the foundation for the gluing of ACyl associative submanifolds in twisted connected sum $G_2$-manifolds presented by the author in \cite{Bera2022}. We study the moduli space of ACyl associative submanifolds with a fixed asymptotic holomorphic curve and a fixed rate, as well as the moduli space where this asymptotic data is allowed to vary, each endowed with a natural topology. We also express their virtual dimensions in terms of data on the asymptotic cross-sections.
\end{abstract}
\section{Introduction} 
 In pursuit of defining enumerative invariants of $G_2$-manifolds that remain unchanged under deformations of the $G_2$-metric, \citet{Joyce2016},  \citet{Doan2017d} have proposed frameworks based on counting closed associative submanifolds, three-dimensional submanifolds calibrated by the defining 3-form of the $G_2$-structure. These are volume-minimizing within their homology class and are thus minimal submanifolds, serving as analogues of holomorphic curves and special Lagrangian submanifolds in Calabi--Yau 3-folds.  \citet[Section~3]{Donaldson1998} suggested an alternative approach to defining invariants via counting $G_2$-instantons. However, such instantons may degenerate by bubbling along associative submanifolds, highlighting the crucial role these submanifolds play in this enumerative theory \cite{Donaldson2009}. Despite the inherent difficulty of constructing $G_2$-manifolds, \citet{Kov03} introduced one of the most powerful methods to date: the twisted connected sum (TCS) construction, which involves gluing asymptotically cylindrical (ACyl) $G_2$-manifolds. As a result, producing associative submanifolds within TCS $G_2$-manifolds via similar gluing methods offers a promising testing ground for these enumerative theories. Such a construction has been explored by the author in \cite{Bera2022}. Part of this article establishes the groundwork for the gluing theorem presented in \cite{Bera2022}. Indeed, the proof of that gluing theorem, along with the necessary conditions for the gluing hypothesis, motivates a detailed study of the deformation theory and moduli space of ACyl associative submanifolds in an ACyl $G_2$-manifold. 

Deformation theory for compact associative submanifolds, with or without boundary, has already been studied by \cites{McLean1998, Gayet2011}, while in this work we study its ACyl version. We investigate the moduli space of ACyl associative submanifolds in an ACyl $G_2$-manifold, focusing first on the case where the asymptotic cross-section and decay rate are fixed, and then on the case where this asymptotic data is allowed to vary. Each moduli space is endowed with a natural topology. Specifically, we consider the moduli space $\mathcal{M}^\mu_{\operatorname{ACyl},\Sigma}$ of ACyl associative submanifolds asymptotic to a holomorphic curve $\Sigma$ with a fixed decay rate $\mu  < 0$ in an ACyl $G_2$-manifold (see \autoref{def ACyl asso}). In \autoref{sec Moduli space of ACyl associative submanifolds}, we endow this moduli space with the $C^\infty_\mu$-topology (see \autoref{def Acyl Cinfty mu topo}) and establish \autoref{thm moduli ACyl asso fixed}, which describes the local structure of the moduli space and computes its virtual dimension i.e., the index of the associated deformation operator, in terms of spectral data of a Dirac operator that encodes the deformation theory of $\Sigma$.

The virtual dimension of the moduli space with fixed asymptotic data, as described above, is always non-positive, and can become negative even when the decay rate is arbitrarily small, if the cokernel of the deformation operator receives contributions from deformations of the asymptotic holomorphic curve. To overcome this issue, we instead allow the asymptotic data to vary and consider the moduli space $\mathcal{M}_{\operatorname{ACyl}}$ of all ACyl associative submanifolds, endowed with the weighted $C^\infty$-topology (see \autoref{def Acyl weighted topo}). In this setting, the cokernel of the corresponding deformation operator no longer acquires contributions from deformations of the asymptotic holomorphic curve, so the virtual dimension becomes non-negative. This is essentially the content of \autoref{thm moduli ACyl asso varyinng} in \autoref{sec Moduli space of ACyl associative submanifolds}, which describes the local structure of this moduli space near ACyl associative submanifolds whose asymptotic holomorphic curves are Morse--Bott i.e., all infinitesimal deformations of the holomorphic curve are integrable (see \autoref{def Morse--Bott}).

The moduli space $\cM^{\hol}_Z$ of holomorphic curves in a Calabi--Yau 3-fold $Z$ is generally not smooth everywhere and is thus best viewed as a complex analytic space. To study the local structure of the moduli space of ACyl associative submanifolds in greater generality than the treatment in \autoref{thm moduli ACyl asso varyinng}, we explain in \autoref{rmk moduli ACyl asso varyinng} that, by equipping $\cM^{\hol}_Z$ with its canonical minimal Whitney stratification, one can decompose $\mathcal{M}_{\operatorname{ACyl}}$ into sub-moduli spaces according to which stratum the asymptotic cross-section lies in, for each of which the local structure can be described as above.

We would like to point out that there are related works addressing deformation problems for ACyl special Lagrangian submanifolds \cite{Salur2010}, coassociative submanifolds \cite{Salur2005}, and Cayley submanifolds \cite{Ohst2015}. All of these, including the present article, share the common feature of following the analysis of elliptic operators on non-compact manifolds developed by \citet{Lockhart1985}. The main difference lies in the specific deformation operators involved, which lead to different type of results concerning unobstructedness. 

\paragraph{Acknowledgements.}I thank my PhD supervisor, Thomas Walpuski, for his constant support and guidance throughout this work, which was earlier part of arXiv:2209.00156v3. I also appreciate Johannes Nordström, Jason Lotay, Dominik Gutwein, Viktor Majewski, and the anonymous referees for their feedback on earlier versions. This material is based upon work supported by the \href{https://sites.duke.edu/scshgap/}{Simons Collaboration on Special Holonomy in Geometry, Analysis, and Physics}. 

\paragraph{Convention.} \textit{Choose} a cut-off function $\chi\in C^\infty(\R, [0,1])$ with $\chi |_{(-\infty,0])}=0$ and $\chi |_{[1,\infty))}=1$.
Set  $\chi_T(t):=\chi(t-T).$
\section{Preliminaries}
In this section we review definitions and basic facts about $G_2$-manifolds and associative submanifolds, which are important to understand this article. To delve further into these topics, we refer the reader to \cite{Joyce2007} and other relevant sources mentioned throughout the discourse.
\subsection{Asymptotically cylindrical (ACyl) $G_2$-manifolds}
A  $3$-form $\phi$ on a $7$-dimensional manifold $Y$ is called \textbf{definite} if the bilinear form $G_\phi: S^2TY\to \Lambda^7(T^*Y)$ defined by 
$G_\phi(u,v):=\iota_u\phi\wedge\iota_v \phi\wedge \phi$
is definite. It uniquely  defines a Riemannian metric $g_\phi$ and a volume form $\vol_{g_\phi}$ on $Y$ satisfying the identity: $G_\phi=6g_\phi \otimes \vol_{g_\phi}$. Moreover it defines 
\begin{itemize}
\item a \textbf{cross product} $\times:\Lambda^2(TY)\to TY$, given by $\phi(u,v,w):=g_\phi(u\times v,w),$	
\item an \textbf{associator} $[\cdot,\cdot,\cdot]:\Lambda^3(TY)\to TY$, given by $[u,v,w]:=(u\times v)\times w+\inp{v}{w}u-\inp{u}{w}v,$
\item a \textbf{$4$-form} $\psi:=*_{g_\phi} \phi\in \Omega^4(Y)$, or equivalently given by $\psi (u,v,w,z):=g_\phi([u,v,w],z)$.\qedhere
\end{itemize} 
 \begin{definition}\label{def G2 manifold}
 A \textbf{$G_2$-manifold} is a $7$-dimensional manifold $Y$ equipped with a torsion-free $G_2$-structure, that is, equipped with a \textbf{definite} $3$-form $\phi\in \Omega^3(Y)$ such that $\nabla_{g_\phi}\phi=0$,
 or equivalently $d\phi=0$ and $ d\psi=0.$
 \end{definition}

\begin{definition}\label{def ACyl G2}
	Let $(Z,\omega,\Omega)$ be a compact Calabi--Yau $3$-fold, where $\omega$ is the K\"ahler form and $\Omega$ is the holomorphic volume form. A $G_2$-manifold $(Y,\phi)$ is called an \textbf{asymptotically cylindrical (ACyl) $G_2$-manifold} with asymptotic cross section $(Z,\omega,\Omega)$ and rate $\nu<0$ if there exist 
	\begin{itemize}
	\item a compact submanifold $K_Y\subset Y$ with boundary 
	\item  a diffeomorphism $\Upsilon:\R^+\times Z\to Y\setminus K_Y$ and  a $2$-form $\varrho$ on $\R^+\times Z$ such that $\Upsilon^*\phi=dt\wedge\omega+\operatorname{Re}\Omega+d\varrho$ with
$$ \abs{\nabla^k\varrho}=O(e^{\nu t})\ \text{as} \ t\to \infty, \forall k\in \N\cup\{0\}.$$ 
	\end{itemize}
 Here $t$ denotes the coordinate on $\R^+$, $\abs{\cdot}$ and Levi-Civita connection $\nabla$ are induced by the product metric on $\R^+\times Z$.
\end{definition}
\begin{remark}Since ACyl $G_2$-manifolds are Ricci flat and complete, by the Cheeger--Gromoll splitting theorem they can have at most two ends, provided they are connected \cite[Theorem 1.1]{Salur2006}. Moreover, if there are two ends, then the holonomy group reduces to a proper subgroup of $G_2$. It was shown in \cite[Theorem 3.8]{Nordstrom2008} that  $\operatorname{Hol}(Y)=G_2$ if and only if the fundamental group $\pi_1(Y)$ is finite and neither $Y$ nor any double cover of $Y$ is homeomorphic to a cylinder. Examples of ACyl $G_2$-manifolds with holonomy exactly $G_2$ are hard to find and the first examples of this kind were found by  \citet[Section 4-5]{KN10}.
\end{remark}

\begin{remark}\label{rmk product ACyl G_2}
	Let $(V,\omega,\Omega)$ be an ACyl Calabi--Yau $3$-fold with asymptotic cross section $S^1\times X$, where $(X,\omega_1,\omega_2,\omega_3)$ is a closed hyperk\"ahler $4$-manifold. Then $$(Y:=S^1\times V,\phi:=d\theta\wedge\omega+\operatorname{Re}\Omega)$$ is an ACyl $G_2$-manifold with asymptotic cross section $$(S^1\times S^1\times X,ds \wedge d\theta+\omega_3, (d\theta-ids)\wedge(\omega_1+i\omega_2)).$$
 Here $\theta$ and $s$ denote the coordinates on the above unit circles. Again by the Cheeger--Gromoll splitting theorem, connected ACyl Calabi--Yau $3$-folds whose holonomy is exactly $\SU(3)$ can have at most one end. 
\end{remark}

\subsection{Asymptotically cylindrical (ACyl) associative submanifolds}\label{subsec Asymptotically cylindrical (ACyl) associative submanifolds}
Every $G_2$-manifold carries a natural calibration, which in turn gives rise to a distinguished class of $3$-dimensional calibrated submanifolds known as associative submanifolds; see \cite[Section 12.3]{Joyce2007}.
\begin{definition}\label{def closed associative}Let $(Y,\phi)$ be a $G_2$-manifold. A $3$-dimensional oriented embedded submanifold $P$ of $Y$ is called an \textbf{associative submanifold} if it is calibrated by the $3$-form $\phi$, that is, $\phi_{|_P}$ is the volume form $\vol_{P,g_\phi}$ on $P$, or equivalently  $\phi_{|_P}$ is the orientation and $[u,v,w]=0$, for all $x\in P$ and $u,v,w\in T_xP$; the later is also equivalent to $T_xP$ being a subalgebra of the cross-product algebra $T_xY$.
\end{definition} 

\begin{example}Let $(Z, \omega, \Omega)$ be a compact Calabi--Yau $3$-fold and $S^1\times Z$ be the product $G_2$-manifold. For any embedded holomorphic curve $\Sigma \subset Z$ and embedded special Lagrangian (i.e. calibrated by $\Re\Omega$) $L\subset Z$ the $3$-dimensional submanifolds $S^1\times \Sigma$ and $\{e^{i\theta}\}\times L$ are associative submanifolds of $S^1\times Z$ for all $e^{i\theta}\in S^1$.
\end{example}	
 
\begin{definition}[Associative cylinders and tubular neighbourhood maps]\label{def asso cylinder}
Let $(Z,J,\omega,\Omega)$ be a compact Calabi--Yau $3$-fold in which $J$ is the complex structure, $\omega$ is the Kähler form and $\Omega$ is the holomorphic volume form. Then 
$$(Y_0:=\R\times Z, \phi_0:=dt\wedge\omega+\operatorname{Re}\Omega)$$
 is a $G_2$-manifold.  
 Let $\Sigma$ be a $2$-dimensional oriented closed submanifold of $Z$ and let $C$ be the cylinder $\R\times \Sigma$. The Calabi--Yau metric on $Z$ induces a metric $g_\Sigma$ on $\Sigma$, the metric $g_{C}=dt^2+g_{\Sigma}$ on $C$, the euclidean normal bundle $NC$ of $C$ with the normal metric $g_{NC}$ and the normal connection $\nabla^\perp_{C}$. Let $\pi:\R\times Z\to Z$ be the standard projection. Then $NC=\pi^*N\Sigma$, the pullback of the normal bundle $N\Sigma$ of $\Sigma$ in $Z$, $g_{NC}=\pi^*g_{N\Sigma}$ and $\nabla^\perp_{C}=\pi^*\nabla^\perp_{\Sigma}.$ Let $\Upsilon_\Sigma:V_\Sigma\subset N\Sigma\to U_\Sigma\subset Z$ be a tubular neighbourhood map of $\Sigma$. 
	Define $$V_C:= \pi^*V_\Sigma\subset NC, \ \ \ U_C:=\R\times U_\Sigma \subset Y_0.$$ 
	The diffeomorphism $\Upsilon_C:V_C\to U_C$ 
	 defined by $$\Upsilon_C((t,\sigma),v):=\big(t,\Upsilon_\Sigma(\sigma,v)\big), \ \ t\in \R, \sigma\in \Sigma, v\in N\Sigma,$$ 
	 is a tubular neighbourhood map of the cylinder $C$. Actually, all the translation invariant tubular neighbourhood maps are of this form.
	 
The cylinder $C=\R\times \Sigma$ is called an \textbf{associative cylinder} in $Y_0=\R\times Z$ if it is an associative submanifold. This is equivalent to requiring that $\Sigma$ be calibrated by the K\"ahler form $\omega$; by Wirtinger's inequality, this is in turn equivalent to $\Sigma$ being a holomorphic curve in $(Z,J)$.
\end{definition}

Now we give the definition of an asymptotically cylindrical associative submanifold.

\begin{definition}\label{def ACyl asso} Let $(Y,\phi)$ be an ACyl $G_2$-manifold with asymptotic cross section $(Z,\omega,\Omega)$ and rate $\nu<0$ as in \autoref{def ACyl G2}. Let $C=\R\times \Sigma$ be a cylinder (possibly disconnected) in $Y_0=\R\times Z$. Let $\Sigma=\amalg_{i=1}^m \Sigma_i$ be the decomposition of $\Sigma$ into connected components, and subsequently $C=\amalg_{i=1}^mC_i$, where $C_i=\R\times \Sigma_i$. Let $\Upsilon_C:V_C\to U_C\subset \R\times Z$ be a translation invariant tubular neighbourhood map of $C$ as in \autoref{def asso cylinder}. 

 A smooth three dimensional oriented submanifold $P$ of $Y$ is said to be an \textbf{asymptotically cylindrical (ACyl) submanifold} with asymptotic cross section $\Sigma$ and rate $\mu=(\mu_1,\mu_2,...,\mu_m)$ with $\nu\leq \mu_i<0$ for all $i=1,2...,m$ if there exist 
 \begin{itemize}
 	\item a compact submanifold with boundary $K_P$ of $P$,
 	\item  a constant $T_0>0$, and a smooth embedding 
	$\Psi_P:(T_0,\infty)\times \Sigma\to U_C\subset \R^+\times Z$
	such that $\Upsilon\circ \Psi_P:(T_0,\infty)\times \Sigma\to Y$ is a diffeomorphism onto $P\setminus K_P$ and 
		$\Psi_P=\Upsilon_C\circ\alpha$ over $(T_0,\infty)\times \Sigma$ for some smooth section $\alpha$ of the normal bundle $NC$ of $C$ which lies in $V_C$ and
		 \begin{equation}\label{eq ACyl asso}
		 	\abs{(\nabla^\perp_{C_i})^k\alpha}=O(e^{\mu_i t}) \  \text{as} \  t\to \infty, i=1,2...,m,\ \forall k\in \N\cup\{0\}.
		 			 	 \end{equation}
		 	  \end{itemize}
 Here $\nabla^\perp_{C}$ is the normal connection on $NC$ induced from the Levi-Civita connection on $\R^+\times Z$ and $\abs{\cdot}$ is taken with respect to the normal metric on $NC$ and cylindrical metric on $C$.
	
	
An associative submanifold is said to be an \textbf{asymptotically cylindrical associative submanifold} with asymptotic cross section $\Sigma$ and rate $\mu$ if it is an ACyl submanifold as above.	
\end{definition}
\begin{remark}\label{rmk ACyl asso} We make a few remarks about the above definition.
\begin{enumerate}[leftmargin=*, labelindent=0pt]
\item An ACyl submanifold with rate $\mu$ is also an ACyl submanifold with any rate $\mu^\prime \geq\mu$. Moreover, it is also an ACyl Riemannian manifold with rate $\mu$ with the metric induced by $g_\phi$ corresponding to the $G_2$-structure $\phi$.
\item If $P$ is an ACyl associative submanifold, then the asymptotic cylinder $C$ is an associative cylinder in $Y_0$. Indeed, the associator on $C_i$
\[ \abs{[\cdot, \cdot,\cdot]_{|C_i}}=O(e^{\mu_i T}),\ \text{as}\  T\to \infty.\]
But $\abs{[\cdot, \cdot,\cdot]_{|C_i}}$ is translation invariant and therefore $[\cdot, \cdot,\cdot]_{|C_i}=0$. \qedhere \end{enumerate}
\end{remark}

\subsection{Normal bundles and canonical isomorphisms}
The following sets up our conventions for the normal bundle of a submanifold, the tubular neighbourhood map, and various canonical isomorphisms, which will be used extensively throughout the article to describe the deformation theories of ACyl associative submanifolds. 
\begin{definition}[Normal bundle]\label{def prelim normal bundle}
	Let  $Y$ be a manifold and $M$ be a submanifold of it. The \textbf{normal bundle} $\pi:NM\to M$ is characterised by the exact sequence of vector bundle morphisms
	\begin{equation}\label{eq normal bundle exact seq}
	0\to TM\to TY_{|M}\to NM\to 0.\qedhere
	\end{equation}
\end{definition}	
\begin{definition}[Tubular neighbourhood map]\label{def prelim tubular neighbourhood map}
 A \textbf{tubular neighbourhood map} of $M$ is a diffeomorphism between an open neighbouhood $V_M$ of the zero section of the normal bundle $NM$ of $M$ that is convex in each fiber and an open neighbourhood $U_M$ (tubular neighbourhood) of $M$ in $Y$,
$$\Upsilon_M:V_M\to U_M$$ that takes the zero section $0$ to $M$ and  
the composition $NM\to 0^*TNM\xrightarrow{d\Upsilon_M} TY_{|M}\to  NM$ is the identity.
\end{definition}

\begin{definition}[Canonical extension of normal vector fields]
\label{def canonical extension normal vf}
	The tangent bundle $TNM$ fits into the exact sequence: 
$0\to \pi^*NM\xrightarrow{i} TNM \xrightarrow{d\pi} \pi^*TM\to 0.$
This induces an \textbf{canonical extension map}, which extends normal vector fields on $M$ to vector fields on $NM$:
\[\widetilde{\bullet}:C^\infty(NM)\hookrightarrow \Vect(NM),\ \ \  u\mapsto \tilde u:=i(\pi^*u).\qedhere\]
\end{definition}
\begin{notation}\label{notation u}
There are instances in this article where it is more appropriate to use the notation $\tilde u$ but for simplicity, we will abuse notation and denote it by $u$.
\end{notation}

\begin{remark}\label{rmk commutator}
Since the canonical extensions defined in \autoref{def canonical extension normal vf} simply extend normal vectors fiberwise as constant vectors, it follows that $[\tilde{u}, \tilde{v}] = 0$ for all $u, v \in C^\infty(NM)$. This observation will be useful later.
\end{remark}

\begin{definition}[Canonical isomorphisms]
For any section $u\in C^\infty(NM)$ we define the \textbf{graph of $u$} by
\[\Gamma_u:=\{\big(x,u(x)\big)\in NM:x\in M\}.	
\]
This is a submanifold of $NM$ and the bundle $\pi^*NM_{|\Gamma_u}$ fits into the split exact sequence
\begin{equation*}
	\begin{tikzcd}
0\arrow{r} & T\Gamma_u \arrow{r} & TNM_{|\Gamma_u}\arrow[l, bend right, swap, "d(u\circ\pi)"]\arrow{r}{\id-d(u\circ\pi)} &\pi^*NM_{|\Gamma_u} \arrow{r} &0.
\end{tikzcd}
\end{equation*}
In particular, this induces a canonical isomorphism $N\Gamma_u\cong \pi^*NM_{|\Gamma_u}$. Moreover, the composition $T\Gamma_u\to TNM_{|\Gamma_u}\xrightarrow{d\pi} \pi^*TM_{|\Gamma_u}$ is an isomorphism. 
Let $\Upsilon_M:V_M\to U_M$ be a tubular neighbourhood map of $M$. We define
\[C^\infty(V_{M}):=\{u\in C^\infty(NM):\Gamma_u\subset V_M\}.\]
Let $u\in C^\infty(V_{M})$. The composition $\Upsilon_M \circ u : M \to M_u$ provides the canonical diffeomorphism identifying $M$ with $M_u$. Under this identification, there is a \textbf{canonical} bundle isomorphism 
\begin{equation}\label{eq prelim canonical prelim isomorphisms}
\Theta^M_u:NM \to NM_u
\end{equation}
 defined by the following commutative diagram of bundle isomorphisms:
\begin{equation*}
	\begin{tikzcd}[column sep=50pt]
 NM \arrow{r}{\Theta^M_u}\arrow{d}{\pi^*} & NM_{u} \\
 \pi^*NM_{|\Gamma_u} & N\Gamma_u\arrow{l}{\id-d(u\circ\pi)}\arrow{u}{d\Upsilon_M}.
\end{tikzcd}\qedhere
\end{equation*}
\end{definition}

\begin{definition}[Normal connection] \label{def prelim normal connection}
The choice of a Riemannian metric $g$ on $Y$ induces a splitting of the exact sequence \autoref{eq normal bundle exact seq}, that is,
 $TY_{|M}=TM\perp NM$
Denote by $\cdot^\parallel$ and $\cdot^\perp$ the projections onto the first and second summands respectively. The Levi-Civita connection $\nabla$ on $TY_{|M}$ decomposes as 
\[\nabla=
\begin{bmatrix}
\nabla^\parallel & -\operatorname{II}^*\\
\operatorname{II} &\nabla^\perp
\end{bmatrix}.
\]
Here $\operatorname{II}\in \Hom(S^2TM,NM) $ is the \textbf{second fundamental form} of M, $\nabla^\parallel$ is the Levi-Civita connection on $M$ and $\nabla^\perp$ is the \textbf{normal connection} on $NM$.
\end{definition}

\section{Moduli space of holomorphic curves in Calabi--Yau $3$-folds}\label{holo curves in Calabi--Yau}  
In this section, we study the local structure of $\cM^{\hol}_Z$ the moduli space of embedded closed holomorphic curves in a Calabi--Yau $3$-fold $Z$ and aim to prove \autoref{thm moduli holo in Z}.

\begin{definition}\label{def moduli holo}
Let $Z$ be a Calabi--Yau $3$-fold and $\mathcal S$ be the set of all $2$-dimensional oriented closed smooth embedded submanifolds. 
We define the \textbf{$C^k$-topology} on the set $\mathcal S$ by specifying a basis which is a collection of all sets of the form $\{\Upsilon_{\Sigma}(\Gamma_u):u\in\mathcal V^k_{\Sigma}\}$, where ${\Sigma\in \mathcal S}$, $\Upsilon_{\Sigma}$ is a tubular neighbourhood map of $\Sigma$  and $\mathcal V^k_{\Sigma}$ is an open set in $C^\infty(V_{\Sigma})$. Here the topology on $C^\infty(V_{\Sigma})$ is induced by the $C^k$-norm on $C^\infty(N{\Sigma})$.
The \textbf{$C^\infty$-topology} on the set $\mathcal S$ is the inverse limit topology of $C^k$-topologies on it, that is, a set is open in the $C^\infty$-topology if and only if it is open in every $C^k$-topology. 

 The \textbf{moduli space}  $\cM^{\hol}_Z$ of embedded closed holomorphic curves in $Z$ is the subset of all submanifolds $\Sigma$ in $\mathcal S$ which are $J$-holomorphic in $Z$. The topologies on $\cM^{\hol}_Z$ are the subspace topologies of the above. 
\end{definition}
 
For $\Sigma\in \cM^{\hol}_Z$ we denote the complex structure on $\Sigma$ by $j$, which is just the restriction of $J$. Let $\Upsilon_\Sigma:V_\Sigma\subset N\Sigma\to U_\Sigma\subset Z$ be a tubular neighbourhood map of $\Sigma$. For $u\in C^\infty(V_\Sigma)$ we denote by $\Sigma_u$ the submanifold $\Upsilon_{\Sigma}(\Gamma_u)$. Then $\Sigma_u$ is $J$-holomorphic if and only if $u$ satisfies the following non-linear Cauchy-Riemann equation
 $$0=\bar{\partial}_Ju:=\frac12(du+\Upsilon_{\Sigma}^*J(u) \circ du\circ j)\in C^\infty(\overline{\Hom}_\C(T\Sigma,u^*TV_\Sigma)).$$  
Here $\overline{\operatorname{Hom}}_\C$ represents the complex anti-linear homomorphisms. The linearization of the nonlinear map $\bar{\partial}_J:C^\infty(V_\Sigma)\to C^\infty(\overline{\operatorname{Hom}}_\C(T\Sigma,u^*TV_\Sigma))$ at the zero section is described in \cite[Proposition 3.1.1]{McDuff2012}. This is the linear map $\mathfrak d_{\Sigma,J}:C^\infty(N\Sigma)\to C^\infty(\overline{\operatorname{Hom}}_\C(T\Sigma,T{V_\Sigma}_{|\Sigma}))$, defined by
$$\mathfrak d_{\Sigma,J}\xi:=\frac12(\nabla_\Sigma\xi+J\circ (\nabla_\Sigma\xi)\circ j),\quad  \xi\in C^\infty(N\Sigma).$$
The tangential component of $\mathfrak d_{\Sigma,J}$ can be discarded for the deformation theory. The following normal component actually controls the deformation theory.
\begin{definition}\label{def normal Cauchy-Riemann operator }For $\Sigma\in \cM^{\hol}_Z$ the \textbf{normal Cauchy-Riemann operator} $\bar \partial ^N_{\nabla,J}:C^\infty(N\Sigma)\to C^\infty(\overline{\operatorname{Hom}}_\C(T\Sigma,N\Sigma))\cong \Omega^{0,1}(\Sigma, N\Sigma)$ is defined by
\begin{equation*}
	\bar \partial ^N_{\nabla,J}\xi:=\frac12(\nabla_\Sigma^\perp\xi+J\circ (\nabla_\Sigma^\perp\xi)\circ j).\qedhere
\end{equation*}
\end{definition}
We now study the local structure of the moduli space of holomorphic curves.
\begin{definition} Let $Y_0$ be as in \autoref{def asso cylinder}.
The \textbf{multiplication map}, $\times_{Z}:TZ\times TZ\to TZ$ is defined by the orthogonal projection of the cross product in $Y_0$ onto $TZ$ or equivalently, for all vector fields $x, y, w$ on $Z$
\begin{equation*}g(x\times_{Z}y, w)=\Re\Omega(x, y, w).\qedhere
\end{equation*}
\end{definition}

\begin{remark}\label{rmk cross product on Z}
Let $\Sigma$ be an oriented smooth surface in $Z$. The following are equivalent (up to orientations): 
(i) $\Sigma$ is $J$-holomorphic in $Z$, (ii) $C=\R\times \Sigma$ is associative in $Y_0$,
(iii) for all $x, y\in T\Sigma$ and $w\in TZ$, $\operatorname{Re}\Omega(x, y, w)=0$,
(iv) for all $x, y\in T\Sigma$, $x\times_{Z}y=0$,
(v) for all $x\in T\Sigma$ and $n\in N\Sigma$, $x\times n\in N\Sigma$.
\end{remark}
\begin{definition}\label{def Dirac Sigma Acyl}
Let $\Sigma$ be an embedded closed $J$-holomorphic curve in $Z$. The map $\gamma_\Sigma:T\Sigma\to \overline{\End}_\C(N\Sigma)$ given by 
$\gamma_\Sigma(f_\Sigma)(v_\Sigma):=f_\Sigma\times v_\Sigma$ is a skew symmetric $J$-anti-linear Clifford multiplication. Moreover the normal bundle $N\Sigma$ together with the metric $g_{N\Sigma}$, Clifford multiplication $\gamma_\Sigma$ and the metric connection $\nabla^\perp:=\nabla^\perp_\Sigma$ is a {Dirac bundle}, that is, $\nabla^\perp\gamma_\Sigma=0$. The {Dirac operator} is given by
\begin{equation}\label{eq: Fueter on Riemann surface}
	\mathbf D_\Sigma:=\gamma_\Sigma\circ\nabla^\perp= \sum_{i=1}^2f_i\times \nabla_{f_i}^\perp
\end{equation}
where $\{f_i\}$ is a local orthonormal oriented frame on $\Sigma$. Note that $\mathbf D_\Sigma$ is $J$-anti-linear. The map $\gamma_\Sigma$ induces a $J$-anti-linear isomorphism 
	$\boldsymbol{\gamma}_\Sigma: \overline{\operatorname{Hom}}_\C(T\Sigma,N\Sigma)\to {N\Sigma}$
	which is defined by
\begin{equation*}	
	\boldsymbol{\gamma}_\Sigma(f^*_\Sigma\tn v_\Sigma)=\gamma_\Sigma(f_\Sigma)(v_\Sigma).
\qedhere \end{equation*}
\end{definition}
\begin{remark}\label{rmk Clifford multi Acyl}
	$\boldsymbol{\gamma}_\Sigma$ fits into the following commutative diagram:
\[	
\begin{tikzcd}
C^\infty(N\Sigma) \arrow{r}{\bar \partial ^N_{\nabla,J}} \arrow[swap]{dr}{\mathbf D_\Sigma} & C^\infty(\overline{\operatorname{Hom}}_\C(T\Sigma,N\Sigma)) \arrow{d}{\boldsymbol{\gamma}_\Sigma} \\
 &C^\infty(N\Sigma).
\end{tikzcd}
 \]
 The above discussion can also be found in \cite[Section 5.3]{CHNP15} for $S^1 \times \Sigma$ instead of ${\R} \times \Sigma$.
  \end{remark}

\begin{definition}\label{def nonlinear map holo curves}
  Let $\Sigma\in \mathcal S$ and let $\Upsilon_\Sigma:V_\Sigma\subset N\Sigma\to U_\Sigma\subset Z$ be a tubular neighbourhood map of $\Sigma$. We define $\mathcal F:C^\infty(V_\Sigma)\to C^\infty(N\Sigma)$ as follows. If $u\in C^\infty(V_\Sigma)$ and $v\in C^\infty(N\Sigma)$,
 \begin{equation*}
 	\inp{\mathcal F(u)}{v}_{L^2}:=\int_{\Gamma_u}\iota_v(\Upsilon_\Sigma^*\Re\Omega).
 \end{equation*}
The notation $v$ in the integrand is the extension vector field of $v\in C^\infty(N\Sigma)$ in the tubular neighbourhood as in \autoref{notation u}.
\end{definition} 

\begin{prop}\label{prop self ad holo lin}For $u\in C^\infty(V_\Sigma)$, we have $\mathcal F(u)=0$ if and only if $\Sigma_u:=\Upsilon_{\Sigma}(\Gamma_u)$ is $J$-holomorphic.  If $\Sigma$ is $J$-holomorphic then the linearization of $\mathcal F$ at zero, $d\mathcal F_0:C^\infty(N\Sigma)\to C^\infty(N\Sigma)$ 
is given by $$d\mathcal F_0=J\mathbf D_\Sigma.$$
This is a first order formally self-adjoint elliptic operator.
\end{prop}
\begin{proof}The first part follows from the \autoref{rmk cross product on Z}. Let $\{f_1,f_2=Jf_1\}$ be a local oriented orthonormal frame for $T\Sigma$. Then for $u, v \in C^\infty(N\Sigma)$, 
$$\frac d{dt}\big|_{{t=0}}\inp{\mathcal F(tu)}{v}_{L^2}=\frac d{dt}\big|_{{t=0}} \int_{\Gamma_{tu}}\iota_v(\Upsilon_\Sigma^*\operatorname{Re}\Omega)
=\int_{\Sigma}\mathcal L_u\iota_v(\Upsilon_\Sigma^*\operatorname{Re}\Omega).$$ 
This is same as $\int_{\Sigma}\iota_v\mathcal L_u(\Upsilon_\Sigma^*\operatorname{Re}\Omega)+{\iota_{[u, v]}(\Upsilon_\Sigma^*\operatorname{Re}\Omega)}$. Following \autoref{notation u} we have $[u, v]=0$ and therefore this is equal to
$$\int_{\Sigma}\iota_v(\Upsilon_\Sigma^*\operatorname{Re}\Omega)(\nabla_{f_1}u, f_2)+\iota_v(\Upsilon_\Sigma^*\operatorname{Re}\Omega)(f_1,\nabla_{f_2}u)+ \int_{\Sigma}{\iota_{\nabla_vu}(\Upsilon_\Sigma^*\operatorname{Re}\Omega)}+{\iota_v \nabla_u(\Upsilon_\Sigma^*\operatorname{Re}\Omega)}.$$ 
As $\Sigma$ is $J$- holomorphic and $\nabla\Omega=0$ the last two integrals vanish with the help of the \autoref{rmk cross product on Z}. Thus the whole expression above becomes  
$$\int_{\Sigma}-\inp{f_2\times \nabla^{\perp}_{f_1}u}{v}+\inp{f_1\times \nabla^{\perp}_{f_2}u}{v}
=\int_{\Sigma}\inp{J(f_1\times \nabla^{\perp}_{f_1}u)}{v}+\inp{J(f_2\times \nabla^{\perp}_{f_2}u)}{v}.$$
The right hand side is $\int_{\Sigma}\inp{J\mathbf D_\Sigma u}{v}$, which proves the second part of the proposition. Since $\mathbf D_\Sigma$ is a $J$-anti-linear (which is proved in  \autoref{prop Fueter Acyl}~\ref{prop Fueter Acyl anti linear}) and self-adjoint Dirac operator, we obtain the last part of the proposition.
\end{proof}

The following theorem describes the local structure of the moduli space $\cM^{\hol}_Z$. Although the statement itself is not new and already appears in \cite[Proof of Theorem 3.1.5]{McDuff2012} and \cite[Section~1, Main Theorem]{Kodaira1962}, the approach adopted here is different. Our method is specific to Calabi–Yau 3-folds, whereas those references treat a more general setting. The novelty of our approach is that it gives a simpler proof in the present context. More precisely, our argument is in the spirit of \cite[Proof of Theorem 3.1.5]{McDuff2012}, but instead of working with $\bar{\partial}_J$ together with the more involved parallel-transport identifications used there, we use the nonlinear map $\mathcal F$ introduced above. This is closely analogous to the nonlinear map defined in \autoref{def nonlinear ACyl asso map} for the deformation theory of associative submanifolds.

\begin{theorem}\label{thm moduli holo in Z}Let $\Sigma$ be a closed holomorphic curve in a Calabi--Yau $3$-fold $(Z,J)$. Then $\cM^{\hol}_Z$ near $\Sigma$ is homeomorphic to the zero set of a smooth map (obstruction map or Kuranishi map) 
$$\ob_\Sigma:\mathcal I_\Sigma \to \coker\bD_\Sigma,$$
where the operator $\bD_\Sigma$ is the Dirac operator defined in \autoref{eq: Fueter on Riemann surface} and $\mathcal I_\Sigma$ is an open neighbourhood of $0$ in $\ker\bD_\Sigma$.
Moreover the index of the operator $\mathbf D_\Sigma$ is zero and there is a $J$-anti-linear isomorphism $\boldsymbol{\gamma}_\Sigma: \overline{\operatorname{Hom}}_\C(T\Sigma,N\Sigma)\to {N\Sigma}$ such that the normal Cauchy-Riemann operator $\bar \partial ^N_{\nabla,J}$ on $\Sigma$ (see \autoref{def normal Cauchy-Riemann operator }) satisfies $$\boldsymbol{\gamma}_\Sigma\circ \bar \partial ^N_{\nabla,J}=\mathbf D_\Sigma.$$
\end{theorem}
\begin{proof} We observed the fact $\boldsymbol{\gamma}_\Sigma\circ \bar \partial ^N_{\nabla,J}=\mathbf D_\Sigma$ in \autoref{rmk Clifford multi Acyl}.
 Extending the nonlinear map $\mathcal F$ to Hölder spaces we get a smooth map
 $$\mathcal F:C^{2,\gamma}(V_\Sigma)\to C^{1,\gamma}(N\Sigma).$$
\autoref{prop self ad holo lin} implies that the linearization of $\mathcal F$ at zero is an elliptic operator and hence Fredholm. By the implicit function theorem applied to $\mathcal F$ we obtain the map $\ob_\Sigma$ as stated in the theorem (see \cite[Proposition 4.2.19]{Donaldson1990}). We only need to prove that, if $u\in C^{2,\gamma}(V_\Sigma)$ with $\mathcal F(u)=0$, then $u\in C^\infty(V_\Sigma)$. To prove this, we observe  
 $$0=\bD_\Sigma(\mathcal F(u))=a(u,\nabla_\Sigma^\perp u) (\nabla_\Sigma^\perp)^2 u+ b(u,\nabla_\Sigma^\perp u).$$
The operator  $\bD_\Sigma\circ \mathcal F$ is a non-linear elliptic operator after choosing the neighbourhood $V_\Sigma$ sufficiently small. Since $a(u,\nabla_\Sigma^\perp u)\in C^{1,\gamma}$ and  $b(u,\nabla_\Sigma^\perp u)\in C^{1,\gamma}$, by Schauder elliptic regularity (see \cite[Theorem 1.4.2]{Joyce2007}) we obtain $u\in C^{3,\gamma}$. By repeating this argument we get higher regularity, which completes the proof of the theorem. 
\end{proof}

\begin{definition}\label{def Morse--Bott}
A closed holomorphic curve $\Sigma$ in a Calabi--Yau $3$-fold $Z$ is said to be \textbf{Morse--Bott} if the obstruction map $\ob_\Sigma$ vanishes near $0$ in $\mathcal I_\Sigma$.
\end{definition}

\begin{remark}Note that $\ker \mathbf D_{\Sigma}$ is symplectic vector space with the symplectic form $\boldsymbol{\omega}$ defined as follows: for $\xi_1,\xi_2 \in \operatorname{ker}\mathbf D_{\Sigma}$, 
\begin{equation}\label{eq symplectic form}
\boldsymbol{\omega}(\xi_1,\xi_2):=\int_\Sigma \omega(\xi_1,\xi_2)\omega.
\end{equation}
If a closed holomorphic curve $\Sigma$ in a Calabi--Yau $3$-fold $Z$ is Morse--Bott  then the moduli space of holomorphic curves near $\Sigma$ is a smooth symplectic manifold of dimension $\dim \ker \mathbf D_\Sigma$. 
\end{remark}

\begin{example}\label{eg special case tori}We consider a special case where the Calabi--Yau $3$-fold $(Z,J,\omega,\Omega)$ is ${T^2}\times X$ with $(X,\omega_1,\omega_2,\omega_3)$ a closed hyperk\"ahler $4$-manifold and
$$\omega=ds \wedge d\theta+\omega_3, \ \ \ \ \Omega=(d\theta-ids)\wedge(\omega_1+i\omega_2),$$
where $(\theta,s)$ denotes the coordinate on $T^2$. For each point $x\in X$, consider the $J$-holomorphic torus $T^2_x:={T^2}\times\{x\}\subset Z$. Then $NT^2_x={T^2}\times T_xX$ with flat normal connection, and hence $ \ker \mathbf D_{T^2_x}\cong T_xX$ is of dimension $4$.
Since already $\{T^2_y:y\in X\} \subset \cM^{\hol}_Z$, $T^2_x$ must be Morse--Bott.
\end{example}

\begin{example}\label{eg special case SL}We again consider the Calabi--Yau $3$-fold $T^2\times X$ as in \autoref{eg special case tori}. For each point $z\in T^2$ and $I_3$-holomorphic curve $\Sigma$ in $X$ consider the holomorphic curve $\Sigma_z:=\{z\}\times\Sigma$ in $Z$. Observe that $\Sigma$ is a special Lagrangian in $(X,\omega_1,\omega_2+i\omega_3)$ with phase $i$. The moduli space of special Lagrangians $\mathcal M_{SL}$ near $\Sigma$ in $X$ is smooth and homeomorphic to an open neighbourhood of $0$ in $\mathcal H^1(\Sigma)$, the space of harmonic $1$-forms on $\Sigma$ \cite[Theorem 3.6]{McLean1998}. 
We have an isomorphism $$N\Sigma_z=T_z{T^2}\oplus N_X\Sigma\cong\R^2\oplus T^*\Sigma,$$ where the last isomorphism $N_X\Sigma \to T^*\Sigma$ is given by $\{\nu\to \iota_{\nu}\omega_1\}$. A direct computation shows that under this isomorphism $\mathbf D_{\Sigma_z}$, $J$ and $J\mathbf D_{\Sigma_z}$ are congruent to  
$$\hat D,\hat J,\hat J \hat D:\Omega^0(\Sigma,\R)\oplus\Omega^0(\Sigma,\R)\oplus\Omega^1(\Sigma,\R)\to\Omega^0(\Sigma,\R)\oplus\Omega^0(\Sigma,\R)\oplus\Omega^1(\Sigma,\R)$$
respectively,
where $$\hat D:=
\begin{bmatrix}
0 & 0 & d^*\\
0 & 0 & *d\\
d & -*d & 0
\end{bmatrix},\ \ 
\hat J:=
\begin{bmatrix}
0 & 1 & 0\\
-1 & 0 & 0\\
0 & 0 & -*
\end{bmatrix},\ \
\hat J \hat D=
\begin{bmatrix}
0 & 0 & *d\\
0 & 0 & -d^*\\
-*d & -d & 0
\end{bmatrix}.
$$
Thus 
$ \ker \mathbf D_{\Sigma_z}\cong\mathcal H^0(\Sigma)\oplus \mathcal H^0(\Sigma) \oplus \mathcal H^1(\Sigma).$
Hence  $\Sigma_z$ is Morse--Bott and the moduli space $\cM^{\hol}_Z$ near $\Sigma_z$ is $\{\Sigma^\prime_{z^\prime}:=\{z^\prime\}\times \Sigma^\prime:(z^\prime,\Sigma^\prime)\in T^2\times U_\Sigma\},$ for some open neighbourhood $U_\Sigma$ of $\Sigma$ in $\mathcal M_{SL}$. \qedhere
	\end{example}

\section{Linear analysis on ACyl associative submanifolds}\label{sec Linear analysis on ACyl associative submanifolds}
We will see in \autoref{sec Moduli space of ACyl associative submanifolds} that the deformation operator (Fueter operator) controlling the deformations of an ACyl associative submanifold is an asymptotically translation-invariant uniformly elliptic operator. Moreover, it is asymptotic to the Fueter operator of the asymptotic associative cylinder. By introducing appropriate weighted function spaces  of normal vector fields on ACyl submanifolds, we obtain in this section the Fredholm theory for such operators, provided the weights avoid the wall of critical weights or rates. The index of the operator changes when the weights cross the wall of critical rates. The index formulas presented in \autoref{thm moduli ACyl asso fixed} and \autoref{thm moduli ACyl asso varyinng} are also established in this section. This theory is standard and originally appeared in \cite{Lockhart1985}. A very good exposition can be found in \cite[Chapter 4]{Marshall2002}.


Let $P$ be an ACyl submanifold asymptotic to an associative cylinder $C=\R\times \Sigma$ with rate $\mu<0$ in an ACyl $G_2$ manifold $(Y,\phi)$ with asymptotic cross section $Z$ and rate $\nu\leq \mu$ as in \autoref{def ACyl asso}. Let $\Sigma=\amalg_{i=1}^m \Sigma_i$ be the decomposition of $\Sigma$ into connected components. Then $C=\amalg_{i=1}^mC_i$, where $C_i=\R\times \Sigma_i$. 

\begin{notation}Let $\lambda,\lambda^\prime\in \R^m$. In the above and later, the notation $\lambda<\lambda^\prime$ means $\lambda_i<\lambda_i^\prime$, for all $i=1,2,..,m$. Similarly we denote $\lambda>\lambda^\prime$ and $\lambda=\lambda^\prime$ by replacing $<$ with $>$ and $=$ respectively. For any $s\in \R$, $\lambda+s:=(\lambda_1+s,\lambda_2+s,...,\lambda_m+s)$. Set $\abs{\lambda}:=\sum_{i=1}^m\abs{\lambda_i}.$
\end{notation}

\begin{definition}[{\cites[Section 5]{McLean1998}[Definition 2.18]{Doan2017d}}]
\label{def D_M}
The \textbf{Fueter operator} $\mathbf D_M:C^\infty(NM)\to C^\infty(NM)$ of an associative submanifold $M$ of a $G_2$-manifold $Y$ is defined by 
	$$\mathbf D_M:=\sum_{i=1}^3e_i\times \nabla_{M,e_i}^\perp$$
	where $\{e_i\}$ is a local oriented orthonormal frame on $M$ and $\nabla_{M}^\perp$ is the normal connection induced by the Levi-Civita connection on $Y$. A direct computation shows that this definition is independent of the choice of local frames. This operator is called a twisted Dirac operator in \cite[Section 5]{McLean1998}.
\end{definition}
For an associative cylinder $C$ it turns out (see \autoref{prop Fueter Acyl}\ref{prop Fueter Acyl decomp}) that $$\mathbf D_C=J\partial_t+\mathbf D_\Sigma.$$

\begin{definition} Let $D_{P}:C^\infty(NP)\to C^\infty(NP)$ be a first order, formally self-adjoint elliptic operator. It is called an \textbf{asymptotically translation invariant uniformly elliptic operator} asymptotic to $\mathbf D_C$ if over $C_\infty:=(T_0,\infty)\times \Sigma$ with the identifications $\Upsilon \circ \Upsilon_C$ from \autoref{def ACyl asso} and the canonical bundle isomorphism defined in \autoref{eq prelim canonical prelim isomorphisms}, we have
$$D_{P}=\sigma(D_{P})\circ \nabla_{C}^\perp+ A_P$$
such that  
$$\abs{\nabla_C^k(\sigma(D_{P})-\sigma(\mathbf D_C))}=O(e^{\mu t})\ \ \text{and} \ \ \ \ \abs{(\nabla_C^\perp)^k(A_P)}=O(e^{\mu t})$$ 
as $t\to \infty$ for all $k\in \N\cup \{0\}$. Here $t$ denotes the coordinate on $\R^+$, and $\abs{\cdot}$ and the normal connection $\nabla_C^\perp$ on $NC$ are induced by the product metric on $\R^+\times Z$. Here $\sigma$ is the symbol of a differential operator.
\end{definition}
 \begin{definition}[{Weighted function spaces}]
 \label{Weighted function spaces}
	For each $\lambda=(\lambda_1,\lambda_2,...,\lambda_m)\in \R^m$, a \textbf{weight function}  
	$w_{P,\lambda}:P\to(0,\infty)$ is a smooth function on $P$ such that for all $x=\Upsilon\circ\Psi_P (t,\sigma)$ with $(t,\sigma)$ in $(T_0,\infty)\times \Sigma_i$,
		$$w_{P,\lambda}(x)= e^{-\lambda_i t}.$$
	Let $\lambda\in\R^m$, $k\geq0$ and  $\gamma\in (0,1)$. For a continuous section $u$ of $NP$ with $k$ continuous derivatives we define the \textbf{weighted Hölder norm} by
	$$ \norm{u}_{C^{k,\gamma}_{P,\lambda}}:=\sum_{j=0}^{k}\norm{w_{P,\lambda}(\nabla_{P}^\perp)^ju}_{L^\infty(NP)}+[w_{P,\lambda}(\nabla_{{P}}^\perp)^ku]_{C^{0,\gamma}{(NP)}}.$$ 
	Here the connection $\nabla_{P}^\perp$ on  $NP$ is the projection of the Levi-Civita connection induced by $g_\phi$ for the decomposition $TY_{|P}=TP\oplus NP$  and $\abs{\cdot}$ is respect to the metric $g_\phi$. The \textbf{weighted Hölder space} $C^{k,\gamma}_{P,\lambda}$ is then the  Banach space of continuous sections of $NP$ with $k$ continuous derivatives and finite {weighted Hölder norm} $\norm{\cdot}_{C^{k,\gamma}_{P,\lambda}}$. 
		
	Furthermore, we define the \textbf{weighted $C^\infty$-space $C^{\infty}_{P,\lambda}$} by
	\begin{equation*}C^{\infty}_{P,\lambda}:=\bigcap_{k=0}^\infty C^{k,\gamma}_{P,\lambda}.\qedhere
	\end{equation*}	
\end{definition}

\begin{prop}[{\cite[Theorem 4.2, Theorem 4.3]{Marshall2002}}]\label{prop Sobolev embedding and cptness acyl}
	Let $\lambda,\lambda_1,\lambda_2 \in \R^m$, $\ k,l\in \N\cup\{0\}$,  $\alpha,\beta\in(0,1).$  If $k+\alpha\geq l+\beta$  then the inclusion ${C^{k,\alpha}_{P,\lambda_1}}\hookrightarrow {C^{l,\beta}_{P,\lambda_2}}$ is a continuous embedding, provided $\lambda_1\leq\lambda_2$.
		Moreover, the embedding is compact, provided the inequality is a strict inequality.
\end{prop}

An asymptotically translation invariant uniformly elliptic operator $D_{P}:C_c^\infty(NP)\to C_c^\infty(NP)$ asymptotic to $\mathbf D_C$ has canonical extensions to weighted function spaces and we denote them as follows:
\begin{equation*} D_{P,\lambda}^{k,\gamma}:C^{k+1,\gamma}_{P,\lambda}\to C^{k,\gamma}_{P,\lambda}.\qedhere\end{equation*}

\begin{prop}[{\cite[Theorem 4.6]{Marshall2002}}] \label{prop elliptic regularity Acyl} Suppose $\lambda\in \R^m$, $k\geq0$, $\gamma\in(0,1)$. Let $u,v \in L^1_{\loc}$ be two locally integrable sections of $NP$ such that $u$ is a weak solution of the equation $D_{P}u=v$. 
If $v\in C^{k,\gamma}_{P,\lambda}$, then $u\in C^{k+1,\gamma}_{P,\lambda}$ is a strong solution and there exists a constant $c>0$ such that 		
		$$\norm{u}_{C^{k+1,\gamma}_{P,\lambda}}\leq c \Big(\norm{D_{P,\lambda}^{k,\gamma}u}_{C^{k,\gamma}_{P,\lambda}}+\norm{w_{P,\lambda} u}_{L^\infty{(NP)}}\Big).$$
\end{prop}

To understand the Fredholm theory for an asymptotically translation invariant uniformly elliptic operator $D_{P}$ asymptotic to $\bD_{C}$, we need a small preparation.

\begin{definition}\label{def critical rates} Let $C=\amalg_{i=1}^mC_i$ be an associative cylinder, where each $C_i=\R\times \Sigma_i$ is connected. 
For $\lambda\in \R^m$, we define the \textbf{homogeneous kernel} of rate $\lambda$ by
	$$V_\lambda:=\bigoplus_{i=1}^m V_{\lambda_i}:=\{e^{\lambda_i t}\nu_{\Sigma_i}: \nu_{\Sigma_i}\in C^\infty(N\Sigma_i), \ \mathbf D_{C_i}(e^{\lambda_i t}\nu_{\Sigma_i})=0\}.$$
Set  $d_{\lambda}:=\sum_{i=1}^md_{\lambda_i}:=\dim V_{\lambda_i}.$ Then the \textbf{wall of critical rates} $\mathcal D_C\subset \R^m$ is defined by
	$$\mathcal D_C:=\{(\lambda_1,\lambda_2,..,\lambda_m)\in \R^m: d_{\lambda_i} \neq 0\ \text{for some}\ i\}.$$ 
 We define the Banach spaces $C^{k,\gamma}_{C,\lambda}$ over $C$ as in \autoref{Weighted function spaces}, replacing $P$ by $C$, $NP$ by $NC$ and the weight function $w_{P,\lambda}$ by $w_{C,\lambda}:C=\amalg _iC_i\to \R$ where $w_{C,\lambda}(t,\sigma_i)= e^{-\lambda_i t}$, $\sigma_i\in \Sigma_i$.  
\end{definition}	
\begin{prop}\label{prop Fueter Acyl}
	Let $\nu_\Sigma\in C^\infty(N\Sigma)$. Then the following hold.
	\begin{enumerate}[label=(\roman*), leftmargin=*]				
	    \item \label{prop Fueter Acyl decomp} $\mathbf D_C=J\partial_t+\mathbf D_\Sigma$.
		\item \label{prop Fueter Acyl homogeneous decomp} $\mathbf D_C(e^{\lambda t}\nu_\Sigma)=e^{\lambda t}\big(\mathbf D_\Sigma \nu_\Sigma+\lambda J\nu_\Sigma \big)$.
		\item \label{prop Fueter Acyl homogeneous sum decomp} $\mathbf D_C(e^{\lambda t}t^j\nu_\Sigma)=e^{\lambda t}t^j\big(\lambda J\nu_\Sigma +\mathbf D_\Sigma \nu_\Sigma\big)+je^{\lambda t}t^{j-1}J\nu_\Sigma$.
		\item \label{prop Fueter Acyl anti linear}$\mathbf D_\Sigma(J\nu_\Sigma)=-J\mathbf D_\Sigma\nu_\Sigma$.
		\item \label{prop Fueter Acyl homogeneous kernel}
 $V_\lambda=\{e^{\lambda t}\nu_\Sigma: \mathbf D_\Sigma \nu_\Sigma=-\lambda J\nu_\Sigma \}= \{e^{\lambda t}\nu_\Sigma: (J\mathbf D_\Sigma) \nu_\Sigma=\lambda \nu_\Sigma \}$.
		\item \label{prop Fueter Acyl homogeneous symmetry}
$JV_{\lambda}=V_{-\lambda}$ and $d_\lambda=d_{-\lambda}$ for all $\lambda\in \R$.
\end{enumerate}
\end{prop}
\begin{proof}
Observe that $\nabla_{C}^\perp=dt\tn {\partial_t}+\nabla_{\Sigma}^\perp$. Now for a local oriented orthonormal frame $\{f_i\}$ on $\Sigma$ and $v\in C^\infty(NC)=C^\infty(\R,N\Sigma)$ we have 
\begin{align*}
\mathbf D_Cv &=J\partial_tv +\sum_{i=1}^2f_i\times \nabla_{\Sigma,f_i} ^\perp v=J\partial_tv+\mathbf D_{\Sigma} v.
\end{align*}
This proves (i) and (ii). In fact, (iii) follows from (i). Indeed
$$\mathbf D_C(e^{\lambda t}t^j\nu_\Sigma)=e^{\lambda t}t^j\big(\lambda J\nu_\Sigma +\mathbf D_\Sigma \nu_\Sigma\big)+je^{\lambda t}t^{j-1}J\nu_\Sigma.$$

To see (iv), we compute
\begin{align*}
\mathbf D_\Sigma(J\nu_\Sigma)=	\sum_{i=1}^2 f_i\times \big(\nabla_{\Sigma,f_i}^\perp({\partial_t}\times\nu_\Sigma)\big)
&=\sum_{i=1}^2 f_i\times ({\partial_t}\times\nabla_{\Sigma,f_i}^\perp\nu_\Sigma)+\sum_{i=1}^2 f_i\times ({\partial_t}{f_i}\times\nu_\Sigma)\\
&=	-\sum_{i=1}^2{\partial_t}\times ( f_i\times\nabla_{\Sigma,f_i}^\perp\nu_\Sigma)=-J\mathbf D_\Sigma(\nu_\Sigma)
\end{align*}
Finally, (v) follows from (ii), and (vi) follows from (iv).
\end{proof}

\begin{prop}\label{prop no sum homogeneous acyl}
Let $m\in \N\cup \{0\}$ and $v(t,\sigma)=\sum_{j=0}^me^{\lambda t}t^j\nu_{\Sigma,j}\in C^\infty(NC)$ with $\nu_{\Sigma,j}\in  C^\infty(N\Sigma)$ and $\nu_{\Sigma,m}\neq 0$. If $\mathbf D_Cv(t,\sigma)=0$, then $m=0$ and thus 
$v(t,\sigma)=e^{\lambda t}\nu_{\Sigma,0}.$ 	
\end{prop}
\begin{proof}
	If $\mathbf D_Cv(t,\sigma)=0$ then by \autoref{prop Fueter Acyl}\ref{prop Fueter Acyl homogeneous sum decomp} and comparing the coefficients of $e^{\lambda t}t^{j-1}$, $j\geq 1$ we can see that
 $$\mathbf D_\Sigma \nu_{\Sigma,m}+\lambda J\nu_{\Sigma,m}=0\ \ \text{and}\ \ jJ\nu_{\Sigma,j}+\mathbf D_\Sigma \nu_{\Sigma,j-1}+\lambda J\nu_{\Sigma,j-1}=0.$$
 Thus we have \begin{align*}
 m\norm{\nu_{\Sigma,m}}^2_{L^2(\Sigma)}=  m\norm{J\nu_{\Sigma,m}}^2_{L^2(\Sigma)}&=-\inp{J\nu_{\Sigma,m}}{\mathbf D_\Sigma \nu_{\Sigma,m-1}+\lambda J\nu_{\Sigma,m-1}}_{L^2(\Sigma)}\\
 &=-\inp{\mathbf D_\Sigma J\nu_{\Sigma,m}+\lambda \nu_{\Sigma,m}}{\nu_{\Sigma,m-1}}_{L^2(\Sigma)}\\
 &=\inp{J(\mathbf D_\Sigma \nu_{\Sigma,m}+\lambda J\nu_{\Sigma,m})}{\nu_{\Sigma,m-1}}_{L^2(\Sigma)}=0.
\end{align*}
The last inequality uses the fact from \autoref{prop Fueter Acyl}\ref{prop Fueter Acyl anti linear}. Since $\nu_{\Sigma,m}\neq 0$, the above identity implies that $m=0$. 
 \end{proof}

\begin{remark}In \cite[Equation 1.11]{Lockhart1985}, the homogeneous kernel is defined to consist of expressions of the form $\sum_{j=0}^m e^{\lambda t} t^j \nu_{\Sigma,j}$. However, \autoref{prop no sum homogeneous acyl} shows that, in our setting, the homogeneous kernel must instead take the simpler form given in \autoref{def critical rates}.
Moreover, the critical rates are in canonical one-to-one correspondence with the eigenvalues of the self-adjoint elliptic operator $J\mathbf D_\Sigma$. In particular, they form a countable discrete subset of $\R$.
\end{remark}

\begin{prop}[{\cites[Lemma 3.1, Proposition 3.4, Section 3.3.1]{Donaldson2002}[Theorem 5.1]{Mazya1978}}]\label{lem fourier series Donaldson Acyl} $\mathbf D_C=J\partial_t+\mathbf D_\Sigma:C^{k+1,\gamma}_{C,\lambda}\to C^{k,\gamma}_{C,\lambda}$ is invertible  if and only if $\lambda\in \R^m\setminus \mathcal D_C$. 
 Moreover, any element $u\in \ker \mathbf D_C$ has an $L^2$-orthogonal decomposition
$$u=\sum_{\lambda\in \mathcal D_C}e^{\lambda t} u_{\Sigma,\lambda}.$$
where $u_{\Sigma,\lambda}$ are $\lambda$-eigensections of $J\mathbf D_\Sigma$.
 \end{prop}
\begin{prop}\label{prop Fredholm Donaldson Acyl}Let $\lambda\in \R^m$, $k\geq0$ and $\gamma\in(0,1)$. Then $D_{P,\lambda}^{k,\gamma}$ is Fredholm for all $\lambda\in \R^m\setminus\mathcal D_C$. Moreover, for all $\lambda\in \R^m$, 
$\operatorname{Ker}D_{P,\lambda}^{k,\gamma}$ is finite dimensional, independent of $k$ and $\gamma$. 
\end{prop}
\begin{proof} The fact that the operator $D_{P,\lambda}^{k,\gamma}$ is Fredholm with weight $\lambda\in \R^m\setminus\mathcal D_C$ follows from \cite[Proposition 3.6, Section 3.3.1]{Donaldson2002} combined with \cite[Theorem 4.2, Theorem 4.3]{Marshall2002}. Independence of $k$ and $\gamma$	follows from \autoref{prop elliptic regularity Acyl} and \autoref{prop Sobolev embedding and cptness acyl}. 
\end{proof}
\begin{lemma}\label{lem main Fredholm acyl}
Let $\lambda\in \R^m$ and $\lambda_1,\lambda_2$ be two elements in $\R^m\setminus\mathcal D_C$ with $\lambda_1<\lambda<\lambda_2$, $\abs{\lambda_2-\lambda_1}\leq \abs{\mu}$, and assume that there are no critical rates between $\lambda_{1},\lambda_{2}$ except possibly $\lambda$. We define the set
	$$\mathcal S_{\lambda_2}:=\{u\in C^{k+1,\gamma}_{P,\lambda_2}: D_{P,\lambda_2}^{k,\gamma}u\in C^{k,\gamma}_{P,\lambda_1}\}.$$
	Define a linear map $e_{P,\lambda}:V_\lambda\to C^{k+1,\gamma}_{P,\lambda_2}$ by $e_{P,\lambda}(w):=\chi_{T_0}w$ (under the identification of the canonical bundle isomorphisms $NP\cong NC$ over the end $(T_0,\infty)\times \Sigma$). Then there exists a unique linear map $\widetilde i_{P,\lambda}:\mathcal S_{\lambda_2}\to V_{\lambda}$ such that for any $u\in \mathcal S_{\lambda_2}$ we have
	$$ u- e_{P,\lambda}\circ \widetilde i_{P,\lambda} u\in  C^{k+1,\gamma}_{P,\lambda_1}.$$
	 Moreover, the following statements hold.
	\begin{enumerate}[label=(\roman*), leftmargin=*]
	\item  $\mathcal S_{\lambda_2}\subset C^{k+1,\gamma}_{P,\lambda}$, $\operatorname{Ker}D_{P,\lambda_2}^{k,\gamma}=\operatorname{Ker}D_{P,\lambda}^{k,\gamma}$, and $\operatorname{Ker}\widetilde i_{P,\lambda}=C^{k+1,\gamma}_{P,\lambda_1}$. The restriction map 
		 \begin{equation}\label{eq asymp limit map}i_{P,\lambda}:=\widetilde i_{P,\lambda}:\operatorname{Ker}D_{P,\lambda_2}^{k,\gamma}\to V_{\lambda}	
		 \end{equation}
 satisfies: $\operatorname{Ker}i_{P,\lambda}=\operatorname{Ker}D_{P,\lambda_1}^{k,\gamma}$.
	\item  
	 $\mathcal S_{\lambda_2}=C^{k+1,\gamma}_{P,\lambda_1}+\im e_{P,\lambda}$ and the restriction of ${\bD}_{P,\lambda_2}^{k,\gamma}$, denoted by  
  \begin{equation}\label{eq tilde L}
\widetilde{\bD}_{P,\lambda_1}^{k,\gamma}:C^{k+1,\gamma}_{P,\lambda_1}+\im e_{P,\lambda}\to C^{k,\gamma}_{P,\lambda_1}
  \end{equation}
 has the property that
$\operatorname{Ker}D_{P,\lambda_2}^{k,\gamma}=\operatorname{Ker}\widetilde{\bD}_{P,\lambda_1}^{k,\gamma},\ \ \operatorname{Coker}D_{P,\lambda_2}^{k,\gamma}\cong\operatorname{Coker}\widetilde{\bD}_{P,\lambda_1}^{k,\gamma}.$
\item  \textbf{Wall crossing formula} \textup{\cite[Theorem 6.2]{Lockhart1985}}: $\operatorname{Ind}(D_{P,\lambda_2})=\operatorname{Ind}(D_{P,\lambda_1})+d_\lambda$.
	\end{enumerate} 
\end{lemma}
\begin{proof}
	Let $u$ be an element of $\mathcal S_{\lambda_2}$. Set $\tilde u:=\chi_{T_0}u\in C^{k+1,\gamma}(NC)$. Since $\mathbf D_C \tilde u\in C_{C,\lambda_1}^{k,\gamma}$, by \autoref{lem fourier series Donaldson Acyl} there exists a unique $v \in C_{C,\lambda_1}^{k+1,\gamma}$ such that $\mathbf D_C (\tilde u-v)=0$. Define $$\widetilde i_{P,\lambda}(u):=e^{\lambda t}(\tilde u-v)_{\Sigma,\lambda}$$ following the decomposition in \autoref{lem fourier series Donaldson Acyl}. Since $\tilde u-v-e^{\lambda t}(\tilde u-v)_{\Sigma,\lambda}\in C_{C,\lambda_1}^\infty $, it proves (i) and $\operatorname{Ker}D_{P,\lambda_2}^{k,\gamma}=\operatorname{Ker}\widetilde{\bD}_{P,\lambda_1}^{k,\gamma}$. Note that $\widetilde{\bD}_{P,\lambda_1}^{k,\gamma}$ is well-defined as $\abs{\lambda-\lambda_1}\leq \abs{\mu}$.  

Denote by $W$ the closure of $\operatorname{im}\,\widetilde{\bD}_{P,\lambda_1}^{k,\gamma}$ in $C^{k,\gamma}_{P,\lambda_2}$. Then the linear maps
\begin{equation*}
\operatorname{Coker}\widetilde{\bD}_{P,\lambda_1}^{k,\gamma}\longrightarrow C^{k,\gamma}_{P,\lambda_2}/W 
\qandq 
C^{k,\gamma}_{P,\lambda_2}/W\longrightarrow \operatorname{Coker}{\bD}_{P,\lambda_2}^{k,\gamma}
\end{equation*}
are isomorphisms.
Indeed, the first map is bijective since the inclusion $C^{k,\gamma}_{P,\lambda_1}\hookrightarrow C^{k,\gamma}_{P,\lambda_2}$ is continuous with dense image, and $\operatorname{Coker}\widetilde{\bD}_{P,\lambda_1}^{k,\gamma}$ is finite-dimensional. For the second map, the above implies that every class $[v]\in C^{k,\gamma}_{P,\lambda_2}/W$ admits a representative $v\in C^{k,\gamma}_{P,\lambda_1}$. If $v=D_{P,\lambda_2}^{k,\gamma}u$, then $u\in \mathcal S_{\lambda_2}$, and hence $[v]=0$. This shows that the second map is injective, while surjectivity is immediate.
Consequently, the natural map
$
\operatorname{Coker}\widetilde{\bD}_{P,\lambda_1}^{k,\gamma}\longrightarrow \operatorname{Coker}{\bD}_{P,\lambda_2}^{k,\gamma}
$
is an isomorphism, completing the proof of~(ii).
	
 Since $\operatorname{Ind}(\widetilde{\bD}_{P,\lambda_1})=\operatorname{Ind}(D_{P,\lambda_1})+d_\lambda$, (iii) follows from (ii).  
	\end{proof}
\begin{prop}\label{prop index acyl}$\operatorname{Ker}D_{P,\lambda}^{k,\gamma}$ is independent of $\lambda$ in each connected component of $\lambda\in\R^m\setminus\mathcal D_C$. Moreover,
		for all $\lambda\in\R^m\setminus\mathcal D_C$, we have
	\begin{enumerate}[label=(\roman*), leftmargin=*]
	\item $\operatorname{Coker}D_{P,\lambda}\cong \operatorname{Ker}D_{P,-\lambda} $,
	\item If $s>0$ such that $0$ is the only possible critical rate in between $-s$ and $s$, then $$\dim \operatorname{Ker}D_{P,0}=\dim \operatorname{Ker}D_{P,s}, \ \ \text{and} \ \ 
	\operatorname{Ind}(D_{P,s})=\dim \operatorname{Ker}D_{P,s}-\dim \operatorname{Ker}D_{P,-s}=\frac{d_{0}}{2},$$
	
	\item $\operatorname{Ind}(D_{P,\lambda})=\displaystyle\sum_{\lambda_i\geq 0}\Big(\frac{d_{0,i}}{2}+\displaystyle\sum_{\zeta_i\in\mathcal D_{C_i}\cap(0,\lambda_i)}d_{\zeta_i}\Big)-\sum_{\lambda_i< 0}\Big(\frac{d_{0,i}}{2}+\sum_{\zeta_i\in\mathcal D_{C_i}\cap(\lambda_i,0)}d_{\zeta_i}\Big)$.
		\end{enumerate}
		\end{prop}
\begin{proof}The statement $\operatorname{Ker}D_{P,\lambda}^{k,\gamma}$ is independent of $\lambda$ in each connected component of $\lambda\in\R^m\setminus\mathcal D_C$ is a direct consequence of part (i) of \autoref{lem main Fredholm acyl}.
Since $\bD_P$ is formally self-adjoint, (i) follows with the help of \cite[(4.28)]{Marshall2002} and part (i) of \autoref{lem main Fredholm acyl}. To see (ii), observe from \autoref{lem main Fredholm acyl} that, $\operatorname{Ker}D_{P,s}=\operatorname{Ker}D_{P,0}$
and  $$\dim \operatorname{Ker}D_{P,s}-\dim \operatorname{Ker}D_{P,-s}=\dim \operatorname{Ker}D_{P,-s}-\dim \operatorname{Ker}D_{P,s}+d_0.$$  
Finally, (iii) follows from (ii) and the wall crossing formula in \autoref{lem main Fredholm acyl}.
\end{proof}

\begin{example}\label{eg tori indicial}
Consider the special case where $Z$ is $T^2\times X$ as in \autoref{eg special case tori}. 
 For the holomorphic tori $T^2_x$ with $x\in X$, $\mathbf D_{T^2_x}^2$ is the Laplacian on the flat torus $T^2=\R^2/{(2\pi \Z)^2}$. Therefore the set of critical rates is
 $\mathcal D=\{\pm\sqrt \lambda:\lambda \ \text{is an eigenvalue of}\ \Delta^0 _{T^2}\}.$
  If $\lambda>0$, the direct sum of homogeneous kernels, $$V_{\sqrt \lambda}\oplus V_{-\sqrt \lambda}=E_{\Delta^0 _{T^2}}^\lambda \oplus E_{\Delta^0 _{T^2}}^\lambda\oplus E_{\Delta^0 _{T^2}}^\lambda\oplus E_{\Delta^0 _{T^2}}^\lambda,$$ where $E_{\Delta^0 _{T^2}}^\lambda$ is the $\lambda$-eigenspace of $\Delta^0 _{T^2}$. Thus
  $d_{\sqrt\lambda}=d_{-\sqrt\lambda}=2 \dim E_{\Delta^0 _{T^2}}^\lambda.$
It follows from \cite[Section 6.3.4, p. 132]{Marshall2002} that $\mathcal D\cap [-\sqrt 2,\sqrt 2]=\{0,\pm \sqrt 2\}$ and $d_{\sqrt 2}=12$. If $\lambda=0$ then $V_0\cong T_xX$.
\end{example}

\begin{example} Consider the special case where $Z$ is $T^2\times X$ as in \autoref{eg special case SL}. 
  For the holomorphic curve $\Sigma_z$ with $z\in T^2$, we have  $\hat D^2$ is $\Delta^0 _{\Sigma}\oplus \Delta^0 _{\Sigma}\oplus \Delta^1 _{\Sigma}$, where $\hat D$ is congruent to $\bD_{\Sigma_z}$. Therefore the set of critical rates is $\mathcal D=\{\pm\sqrt \lambda:\lambda \ \text{is an eigenvalue of}\ \Delta^0 _{\Sigma}\oplus \Delta^1 _{\Sigma}\}.$
  If $\lambda>0$, the direct sum of homogeneous kernels, $$V_{\sqrt \lambda}\oplus V_{-\sqrt \lambda}=E_{\Delta^0 _{\Sigma}}^\lambda\oplus E_{\Delta^0 _{\Sigma}}^\lambda \oplus E_{\Delta^1 _{\Sigma}}^\lambda,$$ where $E_{\Delta^i _{\Sigma}}^\lambda$ is the $\lambda$-eigenspace of $\Delta^i _{\Sigma}$, Laplacian acting on $i$-forms. If $\lambda=0$, 
  \[V_0\cong H^0(\Sigma,\R)\oplus H^0(\Sigma,\R)\oplus  H^1(\Sigma,\R).\qedhere\]
\end{example}


\section{Moduli space of ACyl associative submanifolds}\label{sec Moduli space of ACyl associative submanifolds}
 In this section we will consider the moduli space of ACyl associative submanifolds (possibly with multiple ends) in an ACyl $G_2$-manifold. We aim to prove \autoref{thm moduli ACyl asso fixed} and \autoref{thm moduli ACyl asso varyinng} about the local structure of these moduli spaces. 
 
 Let $(Y,\phi)$ be an ACyl $G_2$-manifold with asymptotic cross section $(Z,\omega,\Omega)$ and rate $\nu<0$ as in \autoref{def ACyl G2}. Recall that we have denoted in \autoref{holo curves in Calabi--Yau} the set of all oriented two dimensional connected, closed embedded submanifolds of $Z$ by $\mathcal S$ and the moduli space of holomorphic curves in $\mathcal S$ by $\cM^{\hol}_Z$.   
\begin{definition}
Let $ \mathcal S_{\operatorname{ACyl}}$ be the set of all ACyl submanifolds in $Y$ as in \autoref{def ACyl asso}. Denote by  $\mathcal S_{\operatorname{ACyl},\Sigma}^\mu$ the space of all ACyl submanifolds with fixed asymptotic cross section $\Sigma$ and rate $\mu$. Note that if $\Sigma=\amalg_{i=1}^m \Sigma_i$ and $\mu\in \R^m$ then any $P\in \mathcal S_{\operatorname{ACyl},\Sigma}^\mu$ has $m$ ends.
The \textbf{moduli space}  $\mathcal M_{\operatorname{ACyl}}$ of ACyl associative submanifolds is defined by 
$$\mathcal M_{\operatorname{ACyl}}:=\{P\in \mathcal S_{\operatorname{ACyl}}:P \ \text{is an ACyl associative submanifold}\}.$$ 
 The \textbf{moduli space}  $\mathcal M_{\operatorname{ACyl},\Sigma}^\mu$  of all ACyl associative submanifolds with fixed asymptotic cross section $\Sigma$ and rate $\mu$ is defined by 
 \[\mathcal M_{\operatorname{ACyl},\Sigma}^\mu:= \mathcal S_{\operatorname{ACyl},\Sigma}^\mu\cap \mathcal M_{\operatorname{ACyl}}.\qedhere\] 
 \end{definition}
We now define topologies on the moduli spaces $\mathcal{M}_{\operatorname{ACyl},\Sigma}^\mu$ and $\mathcal{M}_{\operatorname{ACyl}}$, which we call the $C^\infty_\mu$-topology and the weighted $C^\infty$-topology, respectively. Each topology is specified by constructing a basis of open sets around each ACyl associative submanifold in the corresponding moduli space.
To do this, we choose an `end-cylindrical' (ECyl) submanifold near a given ACyl associative submanifold and endow it with an ECyl tubular neighbourhood map. The basic open sets are then defined as open subsets contained within this tubular neighbourhood. This procedure defines the $C^\infty_\mu$-topology on $\mathcal{M}_{\operatorname{ACyl},\Sigma}^\mu$.
For the weighted $C^\infty$-topology on $\mathcal{M}_{\operatorname{ACyl}}$, by varying the cross sections, we further construct a canonical family of ECyl submanifolds, each equipped with its own canonical ECyl tubular neighbourhood map. The basic open sets in this case are taken from the corresponding family of tubular neighbourhoods.
\begin{figure}[h]
  \centering
  \scalebox{.6}{\includegraphics{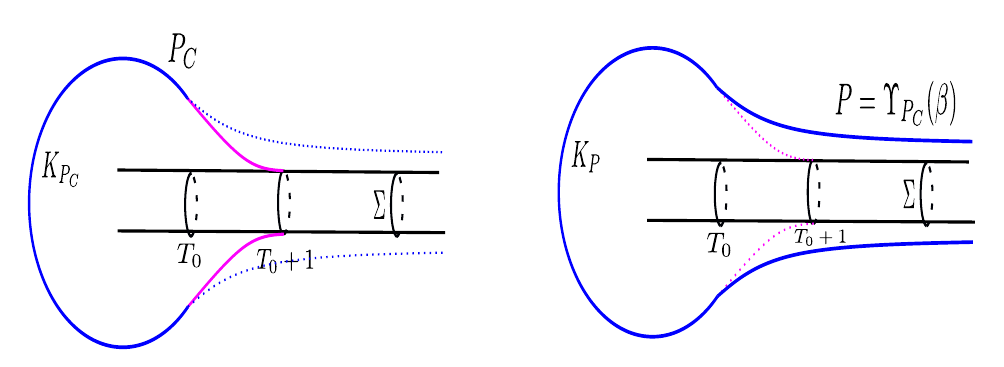}}
  \caption{Construction of ECyl submanifold $P_C$ from an ACyl submanifold $P$.}
  \label{Fig_ACylPC}
\end{figure}

\begin{definition}[ECyl submanifold]\label{def PC Acyl}  Let $P$ be an ACyl submanifold $P$ with asymptotic cylinder $C=\R\times \Sigma$ and $\alpha$ as in \autoref{def ACyl asso}. Recall the diffeomorphism $\Upsilon:\R^+\times Z\to Y\setminus K_Y$ for the end of the ACyl $G_2$-manifold $Y$ and the translation invariant tubular neighbourhood map $\Upsilon_C:V_C\to U_C$ from there. We define an  ECyl submanifold $P_C$ which is diffeomorphic to $P$ but cylindrical on the ends:
$$P_C:=K_P\cup(\Upsilon\circ\Upsilon_C)((1-\chi_{T_0})\alpha).$$  
Set $K_{P_C}:=P_C\setminus \Upsilon((T_0+1,\infty)\times Z).$
\end{definition}

\begin{definition}[ECyl tubular neighbourhood map]\label{def beta}
 Let $P$ be an ACyl submanifold $P$ with asymptotic cross section $\Sigma$ and $\alpha$ as in \autoref{def ACyl asso}. Let $P_C$ be a choice of an ECyl submanifold as in \autoref{def PC Acyl}. A tubular neighbourhood map of $P_C$
	$$\Upsilon_{P_C}:V_{P_C}\to U_{P_C}$$
is said to be end cylindrical (ECyl)  if $V_{P_C}$ and $\Upsilon_{P_C}$ agree with $\Upsilon_*(V_C)$ and $\Upsilon\circ\Upsilon_{C}\circ\Upsilon_*^{-1}$ on $\Upsilon((T_0+1,\infty)\times \Sigma)$ respectively.

Given a choice of an ECyl submanifold $P_C$ and an ECyl tubular neighbourhood map $\Upsilon_{P_C}$, there is a section $\beta$ in $V_{P_C}\subset N{P_C}$ which is zero on $K_P$ and $\Upsilon_*\circ\alpha=\beta\circ \Upsilon$ on $\Upsilon((T_0+1,\infty)\times \Sigma)$ such that $\Upsilon_{P_C}(\Gamma_\beta)$ is $P$.
\end{definition}

\begin{remark}
Set 
$$C_\infty:=(T_0+1,\infty)\times \Sigma,\quad  V_{C_\infty}:=V_{C_{|C_\infty}},\quad  U_{C_\infty}:=U_{C_{|C_\infty}}, \quad P_{C_\infty}:=P_{C_{|\Upsilon(C_\infty)}},\ \  V_{P_{C_\infty}}:=V_{P_{C_{|\Upsilon(C_\infty)}}}.$$ The following commutative diagram helps us to keep track of the definitions above.
\[\begin{tikzcd}
C_\infty\arrow[r,bend left,"\alpha"] \arrow{d}{\Upsilon}&V_{C_\infty} \arrow{l}\arrow{r}{\Upsilon_C} \arrow[swap]{d}{\Upsilon_*} & U_{C_\infty} \arrow{d}{\Upsilon} \\
P_{C_\infty}\arrow[r,bend right,"\beta"]&V_{P_{C_\infty}}=\Upsilon_*(V_{C_\infty})\arrow{l} \arrow{r}{\Upsilon_{P_C}} & Y\setminus K_Y. 
\end{tikzcd}\qedhere
\]
\end{remark}
\begin{definition}[{\textbf{$C^\infty_\mu$-topology}}]\label{def Acyl Cinfty mu topo} 
Given $\mu<0$, we define the \textbf{$C^k_\mu$-topology} on the set $\mathcal{S}_{\operatorname{ACyl},\Sigma}^\mu$ by specifying a basis of open sets around each element $P \in \mathcal{S}_{\operatorname{ACyl},\Sigma}^\mu$. To do so, we choose an end-cylindrical (ECyl) submanifold $P_C$ and an associated ECyl tubular neighbourhood map $\Upsilon_{P_C}$, as in \autoref{def ACyl asso} and \autoref{def PC Acyl}. The basic open sets are of the form
$
\{ \Upsilon_{P_C}(\Gamma_u) : u \in \mathcal{V}^k_{P_C,\mu} \},
$
where $\mathcal{V}^k_{P_C,\mu}$ is an open subset of the weighted function space
$
C^\infty_\mu(V_{P_C}) := C^\infty_{P_C,\mu} \cap C^\infty(V_{P_C}).
$
The topology on $C^\infty_\mu(V_{P_C})$ is induced by the $C^k_\mu$-norm on $C^\infty(NP_C)$.

The \textbf{$C^\infty_\mu$-topology} on $\mathcal{S}_{\operatorname{ACyl},\Sigma}^\mu$ is defined as the inverse limit of the $C^k_\mu$-topologies. That is, a subset is open in the $C^\infty_\mu$-topology if and only if it is open in every $C^k_\mu$-topology.

Finally, the \textbf{$C^\infty_\mu$-topology} on the moduli space $\mathcal{M}_{\operatorname{ACyl},\Sigma}^\mu$ is defined as the subspace topology induced from the $C^\infty_\mu$-topology on $\mathcal{S}_{\operatorname{ACyl},\Sigma}^\mu$.
\end{definition}

\begin{definition}[{\textbf{Weighted $C^\infty$-topology}}]\label{def Acyl weighted topo}
	We define the weighted $C^\infty$-topology on the moduli space $\mathcal M_{\operatorname{ACyl}}$ by specifying a basis of open sets around each element $P \in \mathcal M_{\operatorname{ACyl}}$. To do so, we choose an ECyl submanifold $P_C$ and an associated ECyl tubular neighbourhood map $\Upsilon_{P_C}$, as in \autoref{def ACyl asso} and \autoref{def PC Acyl}. Our aim is to construct a canonical smoothly varying family of ECyl  submanifolds that correspond to variations in the asymptotic cross sections. This is done by constructing a canonical smooth family of diffeomorphisms of $\R\times Z$ as follows.
		
	By \autoref{thm moduli holo in Z}, there is a tubular neighbourhood map $\Upsilon_{\Sigma}:V_{\Sigma}\to U_{\Sigma}$ of ${\Sigma}$ and an obstruction map $\ob_{\Sigma}:\mathcal I_{\Sigma} \to \mathcal O_{\Sigma}$ where $\mathcal I_{\Sigma}:=C^\infty(V_{\Sigma})\cap\operatorname{ker} \mathbf D_{\Sigma}$ and $\mathcal O_{\Sigma}:=\operatorname{coker} \mathbf D_{\Sigma}$. Given an element $ \xi \in \mathcal{I}_{\Sigma} $, we use the extension map $ \widetilde{\bullet}$ from \autoref{def canonical extension normal vf} to obtain a vector field $ \widetilde{\xi} $ on $ V_{\Sigma} $. Applying the differential of the tubular neighbourhood map $ \Upsilon_{\Sigma} $, we obtain the vector field $ d\Upsilon_{\Sigma}(~ \widetilde{\xi}~) $ on $ U_{\Sigma} $. Next, we extend this to a global vector field on $Z $ by multiplying with a cut-off function supported in a neighbourhood of $ U_{\Sigma} $. Finally, we obtain a vector field $ v_{\xi} $ supported on $(T_0, \infty)\times Z\subset \R\times Z$ by translating and multiplying by the cut-off function $ \chi_{T_0} $. In summary, the construction yields a canonical smooth family of vector fields:
	
$$\mathcal I_{\Sigma}\xrightarrow{\widetilde{\bullet}} \Vect (V_\Sigma)\xrightarrow{d\Upsilon_\Sigma} \Vect (U_\Sigma)\hookrightarrow  \Vect (Z)\xrightarrow{\chi_{T_0} \cdot}\Vect (\R\times Z).$$
The time-$1$ flows of this family of vector fields defines a smooth family of diffeomorphisms, $$\Phi: \mathcal I_{\Sigma}\to \Diff(\R\times Z),\ \  \xi\mapsto \Phi_\xi.$$
 In particular, the map $\Phi$ satisfies the following:
\begin{enumerate}[label=(\roman*), leftmargin=*]
\item $\Phi_{0}$ is the identity, and $\Phi_{\xi}$ is the identity on $(-\infty,T_0)\times Z$ for all $\xi\in \mathcal I_{\Sigma}$. 
\item $\Phi_\xi$ is independent of $t$ with $\Phi_\xi(\{t\}\times \Sigma)=\{t\}\times \Upsilon_\Sigma(\xi)$ for all $t\in[T_0+1,\infty)$.
\end{enumerate}

  Given a $\Phi$ as above, we define for any $\xi$ in $\mathcal I_{\Sigma}$, an ECyl submanifold $P^\xi_C$ and an associated ECyl tubular neighbourhood map $\Upsilon_{P^\xi_C}$ :
$$\Upsilon_\xi:=\Upsilon\circ \Phi_\xi,\ \ P^\xi_C:= \big(\Upsilon_{\xi}\circ \Upsilon^{-1}(P_C)\big)\cup K_P$$ 
and
$$\Upsilon_{ P^\xi_C}:=\big(\Upsilon_{\xi}\circ \Upsilon^{-1}\circ \Upsilon_{P_C}\big)\cup \Upsilon_{{P_C}_{|K_P}}:V_{ P_C}\to Y.$$

Finally, the basis of the \textbf{weighted $C^\infty$-topology} on $\mathcal M_{\operatorname{ACyl}}$ is a collection of all subsets of the form 
$$\{\Upsilon_{P^{\xi}_C}(\Gamma_u):u\in\mathcal V_{P_C,\mu},\xi \in \mathcal V_{\Sigma}\}\cap \mathcal M_{\operatorname{ACyl}}$$
 where $\mathcal V_{\Sigma}$ is an open set in $\mathcal I_{\Sigma}$ and $\mathcal V_{P_C,\mu}$ is an open set in $C^\infty_{\mu}(V_{P_C}) := C^\infty_{P_C,\mu} \cap C^\infty(V_{P_C})$.  Here the rate $\mu<0$ is always chosen small enough so that $(\mu,0)\cap \cD_{C_\xi}=\emptyset$ for all associative cylinders $C_\xi:=\R\times \Upsilon_{\Sigma}(\xi)$ with $\xi\in \mathcal I_{\Sigma}$ by choosing $\mathcal I_{\Sigma}$ and $\nu$ small enough if necessary. Anyway this definition will be independent of such choices of $\mu$ (see the last statement of \autoref{thm moduli ACyl asso fixed}). One can also check that it is independent of any other choices that have been made.
 \end{definition}
 
 \begin{figure}[h]
  \centering
  \scalebox{.6}{\includegraphics{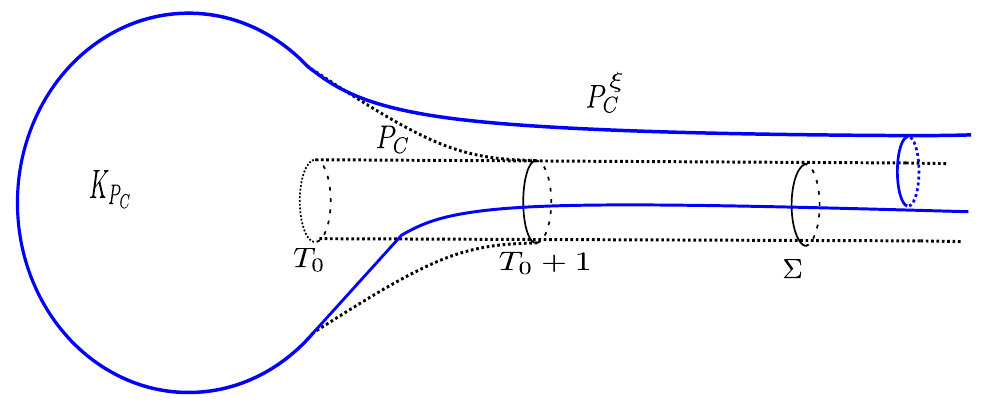}}
  \caption{Construction of $P^\xi_C$, a small ECyl deformation of $P_C$.}
  \label{Fig_PCSigma}
\end{figure}

We define below a nonlinear map whose zero set locally models the moduli space $\mathcal M_{\operatorname{ACyl}}$. Let $P\in \mathcal M_{\operatorname{ACyl}}$ be an ACyl associative submanifold with an asymptotic cylinder $C=\R\times \Sigma$ and rate $\mu$ as in \autoref{def ACyl asso}. We make a choice of an ECyl submanifold $P_C$ and an associated ECyl tubular neighbourhood map $\Upsilon_{P_C}$, as in \autoref{def ACyl asso} and \autoref{def PC Acyl}. There is a \textbf{canonical} bundle isomorphism as defined in \autoref{eq prelim canonical prelim isomorphisms}, which we denote by
\begin{equation}\label{def identification of normal bundle Acyl} 
\Theta_{P}^C:NP_C\to NP.
  \end{equation}
 Given $\xi \in \mathcal I_{\Sigma}$, we also choose ECyl submanifold $P^\xi_C$ and ECyl tubular neighbourhood map $\Upsilon_{P^\xi_C}$ from \autoref{def Acyl weighted topo}.
\begin{definition}\label{def nonlinear ACyl asso map}
We define $\mathfrak F: C^{\infty}_{\mu} (V_{P_C})\times \mathcal I_{\Sigma}\to C^{\infty}{(NP_C)}$ as follows: for all $u\in C^\infty_{\mu} (V_{P_C})$, $\xi \in \mathcal I_{\Sigma}$ and $w\in C_c^\infty (N{P_C})$, 
	  \begin{equation*}
	  	\inp{\mathfrak F(u,\xi)}{w}_{L^2}:=\int_{\Gamma_{u}}\iota_{w}\Upsilon_{P^\xi_C}^*\psi. 
	  	  \end{equation*}
	  	  The notation $w$ in the integrand is actually $\tilde w$, the extension vector field of $w\in C^\infty(NP_C)$ in the tubular neighbourhood as in \autoref{notation u}. The $L^2$ inner product we choose here is coming from the canonical bundle isomorphism $\Theta_{P}^C:NP_C\to NP$ as in \autoref{def identification of normal bundle Acyl} and the induced metric on $NP$ of $g_\phi$.
\end{definition}
\begin{remark}Let $\beta\in C^{\infty}_{\mu} (V_{P_C})$ be the section representing the ACyl associative $P$ as in \autoref{def beta}. Then $\mathfrak F(\beta,0)=0$. Moreover, for $u\in C^{\infty}_{\mu} (V_{P_C})$ and  $\xi \in \mathcal I_{\Sigma}$, we have $\Upsilon_{P^{\xi}_C}(\Gamma_u)\in \mathcal M_{\operatorname{ACyl}}$ if and only if $\mathfrak F(u,\xi)=0$. Hence $\mathcal M_{\operatorname{ACyl}}$ is locally homeomorphic to $\mathfrak F^{-1}(0)$, provided $\mathcal I_{\Sigma}$ and $\mu$ are chosen so that $(\mu,0)\cap \cD_{C_\xi}=\emptyset$ for all associative cylinders $C_\xi:=\R\times \Upsilon_{\Sigma}(\xi)$ where $\xi\in \mathcal I_{\Sigma}$.   
\end{remark}

\begin{prop}\label{prop linearization asso ACyl}For any $\xi\in \mathcal I_{\Sigma}$ the linearization of $\mathfrak F_\xi:=\mathfrak F(\cdot,\xi)$ at $u\in C^{\infty}_{\mu}(V_{P_C})$, 
\begin{equation}\label{eq linearization asso ACyl}
d\mathfrak F_{\xi_{|u}}:C^{\infty}_{P_C,\mu}\to C^{\infty}_{P_C,\mu}
\end{equation}
 is given by the following: for any $v\in C^{\infty}_{P_C,\mu}$ and $w \in C_c^\infty(NP_C)$,
$$\inp{d\mathfrak F_{\xi_{|u}}(v)}{w}_{L^2} =\int_{\Gamma_u}\iota_w\mathcal L_v(\Upsilon_{P^\xi_C}^*\psi)=\int_{\Gamma_u}\biggl\langle\sum_{\substack{\text{cyclic}\\\text{permutations}}}[e_2,e_3,\nabla_{e_1}v],w \biggr\rangle+\int_{\Gamma_u}\iota_{\nabla_wv}(\Upsilon_{P^\xi_C}^*\psi)$$
where $\{e_1,e_2,e_3\}$ is a local oriented orthonormal frame for $T\Gamma_u$ and by abusing notation we are denoting $\Upsilon_{P^\xi_C}^*\nabla$ by $\nabla$ and the associator $\Upsilon_{P^\xi_C}^*[\cdot,\cdot,\cdot]$ by $[\cdot,\cdot,\cdot]$ which are induced from the $3$-form $\Upsilon_{P^\xi_C}^*\phi$. 
\end{prop}
\begin{proof}For the family of sections $\{u+tv:\abs{t}\ll1\}$ we have 
$$\inp{d\mathfrak F_{\xi_{|u}}(v)}{w}_{L^2}=\inp{\frac d{dt}\big|_{{t=0}}\mathcal F_\xi(u+tv)}{w}_{L^2}=\frac d{dt}\big|_{{t=0}} \int_{\Gamma_{u+tv}}\iota_w(\Upsilon_{P^\xi_C}^*\psi).$$
This is equal to $\int_{\Gamma_u}\mathcal L_v\iota_w(\Upsilon_{P^\xi_C}^*\psi)=\int_{\Gamma_u}\iota_w\mathcal L_v(\Upsilon_{P^\xi_C}^*\psi)+{\iota_{[v,w]}(\Upsilon_{P^\xi_C}^*\psi)}$. As $[v,w]=0$, using the definition of Lie derivative this is further equal to 
 \begin{equation*}
 	\int_{\Gamma_u}\sum_{\substack{\text{cyclic}\\\text{permutations}}}\iota_w(\Upsilon_{P^\xi_C}^*\psi)(\nabla_{e_1}v,e_2,e_3)+\int_{\Gamma_u}\iota_{\nabla_wv}(\Upsilon_{P^\xi_C}^*\psi)+\iota_w \nabla_v(\Upsilon_{P^\xi_C}^*\psi).\qedhere
 \end{equation*}
\end{proof}
If the above $u$ represents an associative submanifold then the last integration term vanishes in the above linearization expression and we obtain the following.
\begin{cor}\label{cor linearization}Let $u\in C^{\infty}_{\mu} (V_{P_C})$ and  $\xi \in \mathcal I_{\Sigma}$. If $\Upsilon_{P^{\xi}_C}(\Gamma_u)$ is an ACyl associative submanifold then $d\mathfrak F_{\xi_{|u}}$ is given by the following: for all $v\in C^{\infty}_{P_C,\mu}$ and $w \in C_c^\infty(NP_C)$,
 \begin{equation*}
 	\inp{d\mathfrak F_{\xi_{|u}}(v)}{w}_{L^2} =\int_{\Gamma_u}\biggl\langle\sum_{\substack{\text{cyclic}\\\text{permutations}}}[e_2,e_3,\nabla_{e_1}v],w\biggr\rangle.\qedhere
 	 \end{equation*}
\end{cor}
\begin{prop}\label{prop formal self-adjoint acyl}For any $\xi\in \mathcal I_{\Sigma}$ the linearization $d\mathcal F_{\xi_{|u}}$ of $\mathfrak F_\xi$ at $u\in C^{\infty}_{P,\mu}(V_{P_C})$  is a formally self-adjoint first order differential operator.
\end{prop}
\begin{proof}For all $v, w \in C_c^\infty(NP_C)$, the difference $d\mathcal F_{\xi_{|u}}(v)(w)-d\mathcal F_{\xi_{|u}}(w)(v)$ is
$$\int_{\Gamma_u}\mathcal L_v\iota_w(\Upsilon_{P^\xi_C}^*\psi)-\mathcal L_w\iota_v(\Upsilon_{P^\xi_C}^*\psi)=\int_{\Gamma_u}\iota_w\mathcal L_v(\Upsilon_{P^\xi_C}^*\psi)+{\iota_{[v,w]}(\Upsilon_{P^\xi_C}^*\psi)}-\mathcal L_w\iota_v(\Upsilon_{P^\xi_C}^*\psi).$$
As $[v,w]=0$, this is equal to $\int_{\Gamma_u}\iota_w\iota_v(\Upsilon_{P^\xi_C}^*d\psi)-d(\iota_w\iota_v(\Upsilon_{P^\xi_C}^*\psi))=0$.
\end{proof}

\begin{definition}\label{def L_P acyl}Let $P\in \mathcal M^\mu_{\operatorname{ACyl},{\Sigma}}$ and $\beta\in C^{\infty}_{\mu} (V_{P_C})$ representing $P$ as in \autoref{def beta}. We define $\mathfrak L_P:C_c^\infty(NP_C)\to C_c^\infty(NP_C)$ by
 \[\fL_Pu:= d\mathfrak F_{0_{|\beta}}(u).\qedhere\]
 \end{definition}	

\begin{definition}
	For all $\xi\in \mathcal I_{\Sigma}$ the nonlinear map $Q_\xi: C^{\infty}_{\mu} (V_{P_C})\to C^{\infty}{(NP_C)}$ is defined  by 
	  \begin{equation*}
Q_\xi:=\mathfrak F_\xi-d\mathfrak F_{\xi_{|0}}-\mathfrak F_{\xi}(0).\qedhere 
 \end{equation*}
\end{definition}

\begin{prop}\label{prop ACyl quadratic estimate}There is a constant $c>0$ such that for all $\xi\in \mathcal I_{\Sigma}$, $u, v\in C^{\infty}_{\mu}(V_{P_C})$ and $\eta \in C^{\infty}_{P_C,\mu}$, we have
\begin{enumerate}[label=(\roman*), leftmargin=*]
\item $\abs{d\mathfrak F_{\xi_{|u}}(\eta)-d\mathfrak F_{\xi_{|v}}(\eta)}\leq c (\abs{u-v}+\abs{\nabla^\perp (u-v)})(\abs{\eta}+\abs{\nabla^\perp \eta})$,
\item $\abs{Q_\xi(u)-Q_\xi(v)}\leq c (\abs{u}+\abs{\nabla^\perp u}+\abs{v}+\abs{\nabla^\perp v})(\abs{u-v}+\abs{\nabla^\perp (u-v)})$, 
\item $\norm{Q_\xi(u)}_{C^0_{P_C,2\mu}}\leq c\norm{u}_{C^1_{P_C,\mu}}^2$.
\end{enumerate}
\end{prop}

The proof requires the following lemma whose proof can be found in \cite[Lemma A.1]{Bera2023}.
\begin{lemma}\label{lem quadratic}Let $M$ be a $3$-dimensional oriented submanifold of a $G_2$-manifold $(Y,\phi)$. Let $\Upsilon_M:V_M\subset NM\to Y$ be a tubular neighbourhood map. For all $w \in C^\infty(NM)$ and $u, v, s\in C^\infty(V_M)$  we have\footnote{By abusing notation we denote again by $u,v,w$ the extension over $NM$ of the normal vector fields $u,v,w$ as in \autoref{notation u}. We also denote $\Upsilon_{M}^*\psi$ by $\psi$, where $\psi=*\phi$.}
 $$\abs{\iota_w\mathcal L_u\mathcal L_v\psi}\lesssim \abs{w}\Big(f_{1}\abs{u}\abs{v}
+f_{2}\big(\abs{u}\abs{\nabla v}+\abs{v}\abs{\nabla u}\big)+\abs{\nabla u}\abs{\nabla v}\abs{\psi}\Big),$$
over $\Gamma_s:=\graph s\subset V_M$, where $f_{1}:=\abs{\psi}\abs{\nabla B}+\abs{\psi}\abs{B}^2+\abs{R}\abs{\psi}$ and 
$f_{2}:=\abs{B}\abs{\psi}$.

Here $R$ is the Riemann curvature tensor and $B:C^\infty(NM)\times C^\infty(NM)\to C^\infty(T{V_M})$ is defined by
 $B(u,v):=\nabla_uv$. 
\end{lemma}

 \begin{proof}[{{{{Proof of \autoref{prop ACyl quadratic estimate}}}}}]
For all $w\in C_c^\infty(NP_C)$ and $\eta \in C^{\infty}_{P_C,\mu}$ we write
$$\inp{d\mathfrak F_{\xi_{|u}}(\eta)-d\mathfrak F_{\xi_{|v}}(\eta)}{w}_{L^2}=\int_0^1\frac d{dt}\Big(d\mathcal F_{\xi_{|tu+(1-t)v}}(\eta)(w)\Big) dt$$ and using \autoref{prop linearization asso ACyl} this becomes
$$\int_0^1\Big(\frac d{dt}\int_{\Gamma_{tu+(1-t)v}} \cL_{\eta}\iota_w(\Upsilon_{P^\xi_C}^*\psi)\Big)dt=\int_0^1\int_{\Gamma_{tu+(1-t)v}} \cL_{(u-v)}\cL_{\eta}\iota_w(\Upsilon_{P^\xi_C}^*\psi)dt.$$ 
As $[u-v, w]=0$ and $[\eta,w]=0$ the last expression is same as 
$$\int_0^1\int_{\Gamma_{tu+(1-t)v}} \iota_w \cL_{(u-v)} \cL_{\eta}(\Upsilon_{P^\xi_C}^*\psi)dt.$$
The required estimate in (i) now follows from \autoref{lem quadratic}. 
The estimate in (ii) follows from (i) after writing 
$$Q_\xi(u)-Q_\xi(v)=\int_0^1dQ_{\xi_{|tu+(1-t)v}}(u-v)dt=\int_0^1\big(d\mathfrak F_{\xi_{|tu+(1-t)v}}(u-v)-d\mathfrak F_{\xi_{|0}}(u-v)\big)dt.$$
Finally (iii) follows from (ii). Indeed, substituting $v=0$ we have
\[w_{P_C,2\mu}\abs{Q_\xi(u)}\lesssim w_{P_C,2\mu}(\abs{u}+\abs{\nabla^\perp u})^2\lesssim w^2_{P_C,\mu}(\abs{u}+\abs{\nabla^\perp u})^2.\qedhere\]
\end{proof}
\begin{prop}\label{prop Holder ACyl quadratic estimate}There is a constant $c>0$ such that for all $\xi\in \mathcal I_{\Sigma}$ and $u, v\in C^{k+1,\gamma}_{\mu}(V_{P_C})$ we have
$$\norm{Q_\xi(u)-Q_\xi(v)}_{C^{k,\gamma}_{P_C,2\mu}}\leq c \norm{u-v}_{C^{k+1,\gamma}_{P_C,\mu}}\big(\norm{u}_{C^{k+1,\gamma}_{P_C,\mu}}+\norm{v}_{C^{k+1,\gamma}_{P_C,\mu}}\big).$$
\end{prop}
\begin{proof}
With the above notation and appropriate product operation `$\cdot$', one can express $\mathcal L_u\mathcal L_v\psi$ formally as a quadratic polynomial 
$$\mathcal L_u\mathcal L_v\psi=O(f_1)\cdot u\cdot v+O(f_2)\cdot (u\cdot \nabla^\perp v+v\cdot \nabla^\perp u)+\psi \cdot \nabla^\perp u\cdot \nabla^\perp v$$
with $O(f_{1}):=\psi\cdot \nabla B+\psi\cdot B\cdot B+R\cdot\psi$ and $O(f_2):=B\cdot\psi$, where $B(u,v)=\nabla_uv$. With this observation and similar computations as in the proof of \autoref{lem quadratic} we can prove this Proposition. 
\end{proof}
  
\begin{prop}\label{prop identification of normal bundle Acyl}Let $P\in \mathcal M_{\operatorname{ACyl}}$ be an ACyl associative submanifold and $\beta\in C^{\infty}_{\mu} (V_{P_C})$ representing $P$ as in \autoref{def beta}. Then with the canonical bundle isomorphism $\Theta_{P}^C:NP_C\to NP$ defined in \autoref{def identification of normal bundle Acyl} we have
   \begin{equation*}
\mathfrak L_{P}=(\Theta_{P}^C)^{-1}\circ \mathbf D_{P}\circ \Theta_{P}^C. \qedhere
  \end{equation*}
 \end{prop}
\begin{proof}For all $v,w\in C_c^{\infty}(NP_C)$ we have by definition $\inp{\mathfrak L_{P}v}{w}_{L^2}= \inp{\Theta_{P}^C\mathfrak L_{P}v}{\Theta_{P}^Cw}_{L^2(NP)}$. This is further equal to
	\begin{align*}
 	\int_{\Gamma_\beta}\biggl\langle\sum_{\substack{\text{cyclic}\\\text{permutations}}}[e_2,e_3,\nabla^\perp_{e_1}v],w\biggr\rangle =\int_{P}\biggl\langle\sum_{\substack{\text{cyclic}\\\text{permutations}}}[e_2,e_3,\nabla^\perp_{e_1}(\Theta_{P}^Cv)],\Theta_{P}^Cw\biggr\rangle
 	=\inp{\bD_P \Theta_{P}^Cv}{\Theta_{P}^Cw}_{L^2}.
 	 \end{align*}
 The first equality holds, because $w-\Theta_{P}^Cw\in TP$ and $v-\Theta_{P}^Cv\in TP$, and $[\cdot,\cdot,\cdot]_{|P}=0$.
	 \end{proof}
\begin{prop}\label{prop asymptotically translation elliptic operator}The operators $\mathfrak L_{P}$ and $\mathbf D_P$ are asymptotically translation invariant uniformly elliptic operators and asymptotic to the translation invariant Fueter operator $\mathbf D_C$.  
\end{prop}
\begin{proof}
Since $\mathbf D_P$ is elliptic and by \autoref{prop identification of normal bundle Acyl} $\mathfrak L_{P}=(\Theta_{P}^C)^{-1}\circ\mathbf D_{P}\circ\Theta^C_{P}$,  we obtain that $\mathfrak L_{P}$ is elliptic.

It remains to prove that $\mathfrak L_{P}$ is asymptotic to $\mathbf D_C$. For all $v,w\in C^\infty_c(NC_t)$ where $C_t:=(T_0+1,t)\times\Sigma$ we have
$$\inp{d\mathfrak F_{0_{|0}}\Upsilon_*v}{\Upsilon_*w}_{L^2}=\int_{C_t}\iota_w\mathcal L_v(\Upsilon_{C}^*\Upsilon^*\psi)=\int_{C_t}\iota_w\mathcal L_v(\Upsilon_{C}^*\psi_0)+O(e^{\nu t})\norm{v}_{W^{1,2}(NC)}\norm{w}_{L^2(NC)}.$$
Here $\nu<0$ is the asymptotic rate of the ACyl $G_2$-manifold. By \autoref{prop ACyl quadratic estimate} (i) we obtain $$\inp{\Upsilon^*\mathfrak L_{P}v}{w}_{L^2(NC)}= \inp{\mathbf D_Cv}{w}_{L^2(NC)}+ O(e^{\mu t})\norm{v}_{W^{1,2}(NC)}\norm{w}_{L^2(NC)}.$$
Here $\mu\geq \nu$ is the asymptotic rate of $P$. This completes the proof of the proposition for $\mathfrak L_{P}$ and by \autoref{prop identification of normal bundle Acyl}, the claim for $\mathbf D_P$.
\end{proof}
\begin{definition}\label{def Fueter on ACyl asso} Let $P$ be an ACyl associative submanifold and $\mathbf D_{P}:C^\infty_c(NP)\to C^\infty_c(NP)$ be the Fueter operator defined in \autoref{def D_M}. \autoref{prop asymptotically translation elliptic operator} implies that it is an asymptotically translation invariant uniformly elliptic operator. Therefore it has the following canonical extensions: 
\begin{equation*}
\mathbf D_{P,\lambda}^{k,\gamma}:C^{k+1,\gamma}_{P,\lambda}\to C^{k,\gamma}_{P,\lambda}.
\end{equation*}
Similarly we have canonical extensions:
 \begin{equation*}
 \mathfrak L_{P,\lambda}^{k,\gamma}:C^{k+1,\gamma}_{P_C,\lambda}\to C^{k,\gamma}_{P_C,\lambda}.\qedhere  \end{equation*}
\end{definition}

 We  prove the following theorem about the local structure of the moduli space of ACyl associative submanifolds with a fixed asymptotic cross section and a fixed rate.
 \begin{theorem}\label{thm moduli ACyl asso fixed}Let $P\in \mathcal{M}^\mu_{\operatorname{ACyl},\Sigma}$ be an ACyl associative submanifold. Assume that $\mu$ lies outside the wall of critical rates (see \autoref{def critical rates}). Then  $\mathcal{M}^\mu_{\operatorname{ACyl},\Sigma}$ near $P$ is homeomorphic to the zero set of a smooth map (obstruction map or Kuranishi map) $$\ob_{P,\mu}:\mathcal I_{P,\mu}\to \operatorname{coker}{\mathbf D}_{P,\mu},$$
  where $\mathcal I_{P,\mu}$ is a neighbourhood of $0$ in $\operatorname{ker}{\mathbf D}_{P,\mu}$ and ${\mathbf D}_{P,\mu}$ is the Fueter operator defined in \autoref{def Fueter on ACyl asso}. Moreover, 
  $$\ind \mathbf D_{P,\mu}=-\sum_{i=1}^m\frac{d_{0,i}}{2}-\sum_{i=1}^m\sum_{\lambda_i\in(\mu_i,0)}d_{\lambda_i}\leq 0,$$
   where $d_{\lambda_i}$ is the dimension of the $\lambda_i$-eigenspace of $J\mathbf D_{\Sigma_i} $ and $d_{0,i}:=\dim \ker J\mathbf D_{\Sigma_i}$ (see \autoref{def critical rates}). Moreover $\mathcal M^\mu_{\operatorname{ACyl},\Sigma}= \mathcal M^{\mu^\prime}_{\operatorname{ACyl},\Sigma}$, if $\mu$ and $\mu^\prime$ are lying in the same connected component of the complement of the wall of critical rates.
    \end{theorem}
  
 \begin{proof}
 Extending the nonlinear map $\mathcal F_0$ to weighted Hölder spaces and applying \autoref{prop ACyl quadratic estimate}, \autoref{prop Holder ACyl quadratic estimate} and \autoref{prop linearization asso ACyl} we obtain a well defined map
 $$\mathfrak F_0:C^{2,\gamma}_{\mu}(V_{P_C})\to C^{1,\gamma}_{P_C,\mu}$$
 since, $\mathfrak F_0={\mathfrak F_0}(0)+{d\mathfrak F_0}_{|0}+Q_0$ and $\mathfrak F_0(0)\in C^{1,\gamma}_{P_C,\nu}$ where $\nu<0$ is the asymptotic rate of the ACyl $G_2$-manifold. As $\mu\in (\nu,0)\setminus\mathcal D_C$, by \autoref{prop Fredholm Donaldson Acyl} we have that the linearization of $\mathfrak F_0$ at $\beta$, $\mathfrak L_{P,\mu}$ is a Fredholm operator. Applying the implicit function theorem to $\mathfrak F_0$ and \autoref{prop elliptic regularity Acyl} of elliptic regularity to $\mathfrak L_{P,\mu}$ we obtain the map $\ob_{P,\mu}$ as stated in the theorem (see \cite[Proposition 4.2.19]{Donaldson1990}). We need to only prove that, if $u\in C_{\mu}^{2,\gamma}(V_{P_C})$ with $\mathcal F_0(u)=0$ then $u\in C_{\mu}^{\infty}(V_{P_C})$. To prove this, we observe  
 $$0=\fL_P(\mathcal F_0(u))=a(u,\nabla_{P_C}^\perp u) (\nabla_{P_C}^\perp)^2 u+ b(u,\nabla_{P_C}^\perp u).$$
 Since $a(u,\nabla_{P_C}^\perp u)\in C_{P_C,\mu}^{1,\gamma}$ and  $b(u,\nabla_{P_C}^\perp u)\in C_{P_C,\mu}^{1,\gamma}$, by Schauder elliptic regularity (see \cite[Theorem 1.4.2]{Joyce2007}) and \autoref{prop elliptic regularity Acyl} we obtain $u\in C_{P_C,\mu}^{3,\gamma}$. By repeating this argument we can get higher regularity. The index formula follows from \autoref{prop index acyl}.
  
Finally we prove the last statement of the theorem that the moduli space is the same for two rates $\mu$ and $\mu^\prime$ in the same connected component of the complement of the wall of critical rates. If $\mathfrak F_0(u)=0$, then ${d\mathfrak F_0}_{|0}(u)=-Q_0(u)-{\mathfrak F_0}(0)$. By \autoref{prop ACyl quadratic estimate} we see that ${d\mathfrak F_0}_{|0}$ is also an ACyl uniformly elliptic operator asymptotic to $\mathbf D_C$. Also we have $Q_0(u)\in C^{k,\gamma}_{P_C,2\mu}$ and ${\mathfrak F_0}(0)\in {C^{k,\gamma}_{P_C,\nu}}$. Therefore by an inductive procedure with the help of \autoref{lem main Fredholm acyl} we can conclude that $u\in C^{\infty}_{\mu^\prime} (V_{P_C})$.
\end{proof}

In the theorem above, we observe that the index of the deformation operator is non-positive, indicating that the deformation theory with a fixed cross section is likely obstructed. To obtain a less obstructed deformation theory, we allow the cross sections to vary, which is the subject of the following discussion. Recall that ACyl submanifolds around $P$ with varying asymptotic cross sections are parametrised by $(u,\xi)$ with $u \in C^\infty_{\mu}(V_{P_C})$ and $\xi \in \mathcal I_{\Sigma}$, and that we have a nonlinear map $\mathfrak F$ defined on them in \autoref{def nonlinear ACyl asso map}. In \autoref{prop linearization asso ACyl}, its linearization was computed by fixing the asymptotic cross section parametrised by $\xi$. The following lemma computes the linearization of the nonlinear map when the asymptotic cross sections are allowed to vary, while $\beta$ parametrising $P$ is fixed.

\begin{lemma}\label{lem linearization 2nd component acyl}Let $P$ be an ACyl associative submanifold and $\beta\in C^{\infty}_{\mu} (V_{P_C})$ representing $P$ as in \autoref{def beta}. The linearization of $\mathfrak F(\beta,\cdot)$ at $0\in \mathcal I_{\Sigma}$, $L:\ker \mathbf D_{\Sigma}\to C^{\infty}_{P_C,\mu}$ is given  by
$$L(\xi)=\mathfrak L_P(\chi_{T_0}\Upsilon_*\xi).$$
\end{lemma}
\begin{proof}
For all $w\in C^\infty_c(NP_C)$ and $\xi\in \ker \mathbf D_{\Sigma}$,	
$$\inp{L( \xi)}{w}_{L^2}=\int_{\Gamma_{\beta}}\frac d{dt}\big|_{{t=0}} \iota_w(\Upsilon_{P_C^{t\xi}}^*\psi)=\int_{\Gamma_{\beta}}\iota_w\mathcal L_{(\chi_{T_0}\Upsilon_*\xi)}(\Upsilon_{P_C}^*\psi).$$
Therefore the lemma now follows by a similar line of argument as in \autoref{prop linearization asso ACyl}. \end{proof}
 \begin{definition}\label{def tilde Fueter on ACyl asso}Let $P$ be an ACyl associative and $\Theta^C_P:NP_C\to NP$ be the canonical bundle isomorphism in \autoref{def identification of normal bundle Acyl}. We define 
 $$\widetilde \fL_{P,\mu}:C^\infty_{P_C,\mu}\oplus \ker D_{\Sigma}\to C^\infty_{P_C,\mu}$$
 by $\widetilde \fL_{P,\mu}(u,\xi):=\fL_{P}u+L\xi$, where $\fL_P$ and $L$ are as in \autoref{def L_P acyl} and \autoref{lem linearization 2nd component acyl}. The operator $$\widetilde \bD_{P,\mu}:C^\infty_{P,\mu}\oplus \ker D_{\Sigma}\to C^\infty_{P,\mu}$$ is the operator $\widetilde \fL_{P,\mu}$ under the identification map $\Theta^C_P$, that is,
 \begin{equation*}
 \widetilde \bD_{P,\mu}(u,\xi):= \Theta^C_P\widetilde \fL_{P,\mu}\big((\Theta^C_P)^{-1} u,\xi\big)=\bD_Pu+ \Theta^C_P L\xi.\qedhere
\end{equation*}
\end{definition}
We prove the following theorem about the local structure of the moduli space of ACyl associative submanifolds around an ACyl associative submanifold whose asymptotic cross section is Morse--Bott.
   \begin{theorem}\label{thm moduli ACyl asso varyinng}Let $P$ be an ACyl associative submanifold in an ACyl $G_2$-manifold $Y$ with asymptotic cross section $\Sigma:=\amalg_{i=1}^m\Sigma_i$ is Morse--Bott (see \autoref{def Morse--Bott}). Then the moduli space $\mathcal M_{\operatorname{ACyl}}$ of ACyl associative submanifolds near $P$ is homeomorphic to the zero set of a smooth map (obstruction map or Kuranishi map) $$\widetilde{\ob}_{P,\mu}:\tilde{\mathcal I}_{P,\mu}\to \operatorname{coker}{\widetilde{\mathbf D}}_{P,\mu}\cong\operatorname{ker}{\mathbf D}_{P,\mu}$$
   for some small rate $\mu<0$. Here $\tilde{\mathcal I}_{P,\mu}$ is a neighbourhood of $(0,0)$ in $\operatorname{ker}{\widetilde{\mathbf D}}_{P,\mu}$ and $\widetilde{\mathbf D}_{P,\mu}$ is defined in \autoref{def tilde Fueter on ACyl asso}. Moreover, $$\ind \widetilde{\mathbf D}_{P,\mu}=\frac 12\sum_{i=1}^m\dim \ker \mathbf D_{\Sigma_i}=\sum_{i=1}^m\frac{d_{0,i}}2\geq 0.$$
\end{theorem}
\begin{proof} Suppose $\mathcal I_{\Sigma}$ and $\mu$ are chosen sufficiently small so that $(\mu,0)\cap \cD_{C_\xi}=\emptyset$ for all associative cylinders $C_\xi:=\R\times \Upsilon_{\Sigma}(\xi)$ where $\xi\in \mathcal I_{\Sigma}$. We fix such data. Then near $P$, the moduli space $\mathcal M_{\operatorname{ACyl}}$ is homeomorphic to $\mathfrak F^{-1}(0)$. Since $\Sigma$ is Morse--Bott, for all  $\xi\in \mathcal I_{\Sigma}$ (sufficiently small), the cylinder  $C_\xi$ is associative and therefore $\mathcal F(0,\xi)\in C^{\infty}_{P_C,\mu}$. Using \autoref{prop ACyl quadratic estimate} and \autoref{prop Holder ACyl quadratic estimate} we obtain that the non-linear map $$\mathcal F:C^{2,\gamma}_{\mu}(V_{P_C})\times \mathcal I_{\Sigma} \to C^{1,\gamma}_{P_C,\mu}$$ 
is well defined. By \autoref{lem linearization 2nd component acyl}, the linearization of $\mathcal F$ at $(\beta,0)$, 
$${d\mathfrak F}_{|(\beta,0)}(u,\xi)=\mathfrak L_{P,\mu}u+ L\xi.$$
 Thus ${d\mathfrak F}_{|(\beta,0)}=\widetilde{\mathfrak L}_{P,\mu}$ which is defined in \autoref{def tilde Fueter on ACyl asso}. Now the rest of the proof is similar to the proof of \autoref{thm moduli ACyl asso fixed}. The index formula follows from the proof of \autoref{prop index acyl}.
\end{proof}

\begin{remark}\label{rmk moduli ACyl asso varyinng}
The moduli space of holomorphic curves $\cM^{\hol}_Z$ in the Calabi--Yau $3$-fold $Z$ is generally not smooth everywhere and is therefore best regarded as a complex analytic space. To analyze the local structure of the moduli space of ACyl associative submanifolds in greater generality than in \autoref{thm moduli ACyl asso varyinng}, one can choose the canonical minimal Whitney stratification:
\[
\cM^{\hol}_Z = \bigsqcup_{k \in I} \mathcal{E}^{(k)},
\]
where $\mathcal{E}$ is a smooth manifold of dimension $k$, in the sense of \cite[Chapter I, Sections 1–2]{Gibson1976}. Given such a stratification, the local structure of the sub-moduli space $\mathcal{M}_{\operatorname{ACyl}, \mathcal{E}}$, consisting of ACyl associative submanifolds with asymptotic cross-section $\Sigma$ lying in a fixed stratum $\mathcal{E}$ (assuming a single end for simplicity), can be described following the same reasoning as in the proof of \autoref{thm moduli ACyl asso varyinng}. Specifically, one replaces $\mathcal{I}_{\Sigma}$ with $\mathcal{I}_{\Sigma} \cap (\Upsilon_{\Sigma}^{-1} \mathcal{E})$, and replaces the deformation operator $\widetilde{\mathbf{D}}_{P, \mu}$ with its restriction on $C^\infty_{P_C, \mu} \oplus T_\Sigma \mathcal{E}$, say $\widetilde{\mathbf{D}}_{P, \mu, \mathcal{E}}$. In this case, the index is given by
\[
\operatorname{ind} \widetilde{\mathbf{D}}_{P, \mu, \mathcal{E}} = \dim \mathcal{E} - \frac{1}{2} \dim \ker \mathbf{D}_\Sigma. \qedhere
\]
\end{remark}

\begin{remark}\label{rmk Atiyah-Floer}
Suppose a pair of ACyl associative submanifolds $P_\pm$ satisfies the gluing hypothesis in \cite{Bera2022}, and the asymptotic cross sections $\Sigma_\pm$ of $P_\pm$ are Morse--Bott. Then the moduli spaces $\cM^{\hol}_{Z_\pm}$ of holomorphic curves in $Z_\pm$ are smooth around $\Sigma_\pm$ and are modelled on $\mathcal I_{\Sigma_\pm}$, a neighbourhood of $0$ in $\ker \bD_{\Sigma_\pm}$. The asymptotic limit maps being injective is equivalent to saying $\ker \bD_{P_\pm,\mu} = 0$ for sufficiently small $\mu < 0$ (cf.\ \autoref{lem main Fredholm acyl}), and therefore, by \autoref{thm moduli ACyl asso varyinng} and \autoref{lem main Fredholm acyl}, the moduli spaces $\mathcal M_{\operatorname{ACyl}}(Y_\pm)$ of ACyl associative submanifolds of $Y_\pm$ are smooth around $P_\pm$ and are modelled on ${\mathcal I}_{P_\pm,0}$, a neighbourhood of $0$ in $\operatorname{ker} \mathbf D_{P_\pm,0}$. The canonical maps $\mathcal M_{\operatorname{ACyl}}(Y_\pm) \to \cM^{\hol}_{Z_\pm}$ around $P_\pm$, sending to associated asymptotic cross sections, are therefore locally modelled on a pair of maps ${\mathcal I}_{P_\pm,0} \to \mathcal I_{\Sigma_\pm}$, which are Lagrangian immersions with respect to the symplectic form $\boldsymbol{\omega}$ defined in \autoref{eq symplectic form}. Indeed, one can expect that the linearization of these maps at $0$ are the asymptotic limit maps $ \iota_{P_\pm,0}$ defined in \autoref{eq asymp limit map} of \autoref{lem main Fredholm acyl}, which are already injective here. Set $P := P_\pm$ and $\Sigma := \Sigma_\pm$ . Using an integration by parts argument as in the proof of \autoref{prop formal self-adjoint acyl}, one can see that, if $\xi_k = i_{P,0}(u_k)$ with $u_k \in \operatorname{ker} \bD_{P,0}$, $k = 1,2$, then
\begin{equation*}
\int_{\Sigma} \omega(\xi_1,\xi_2)\,\omega
= \int_{\Sigma} \iota_{\xi_1}\iota_{\xi_2}\!\left(\tfrac{\omega^2}{2}\right)
= \lim_{T\to\infty}\int_{\{T\}\times \Sigma} \iota_{u_1}\iota_{u_2}\Upsilon^*_{P_C}\psi
= \int_{P} \langle \bD_{P}u_1, u_2 \rangle - \langle u_1, \bD_{P}u_2 \rangle
= 0.
\end{equation*}
Therefore, $\operatorname{im} \iota_{P_\pm,0}$ are isotropic. By \autoref{prop index acyl}, $\dim \ker \mathbf D_{P_\pm,0} = \frac{1}{2} \dim \ker \mathbf D_{\Sigma_\pm}$, and this proves that $\operatorname{im} \iota_{P_\pm,0}$ are Lagrangian subspaces.

The gluing hypothesis in \cite{Bera2022} is therefore a Lagrangian intersection hypothesis, which yields a unique closed rigid associative submanifold in the twisted connected sum $G_2$-manifold by the gluing theorem in \cite{Bera2022}. It possibly gives a hint for formulating an Atiyah--Floer type conjecture \cite{Atiyah1988} by considering twisted connected sum $G_2$-manifolds as an analogue of Heegaard splittings of $3$-manifolds. The possible analogy of the Atiyah--Floer conjecture would be some relationship between a yet-to-be-defined Floer homology of closed associative submanifolds, say $HA(Y)$, and Lagrangian Floer homology, say $HF(\mathcal M_{\operatorname{ACyl}}(Y_+), \mathcal M_{\operatorname{ACyl}}(Y_-))$, where $Y = Y_+ \bigcup_Z Y_-$ is the twisted connected sum $G_2$-manifold. A similar analogy with $G_2$-instantons has been pointed out in \cite[Remark 1.7]{SaEarp2013}.
\end{remark}

\printreferences

\end{document}